\documentclass[twoside,leqno,twocolumn]{article}

\usepackage[letterpaper]{geometry}
\usepackage{ltexpprt}

\usepackage{amsmath,amssymb,amsfonts}
\usepackage{algpseudocode}
\usepackage{listings}
\usepackage{booktabs}
\usepackage{graphicx}
\usepackage{tikz}
\usetikzlibrary{arrows.meta}
\usepackage{pgfplots}
\pgfplotsset{compat=1.18}
\usepgfplotslibrary{groupplots}
\usepackage{needspace}
\usepackage{xcolor}
\usepackage[hidelinks]{hyperref}
\usepackage{comment}
\usepackage{xspace}

\definecolor{todocolor}{rgb}{0.75,0.1,0.1}

\newcommand{\R}{\mathbb{R}}
\newcommand{\dd}{\mathop{}\!\mathrm{d}}

\newcommand{\vx}{\boldsymbol{x}}
\newcommand{\vX}{\boldsymbol{X}}
\newcommand{\vtheta}{\boldsymbol{\theta}}

\newcommand{\dfem}{\ensuremath{\partial\mathrm{FEM}}\xspace}

\def\lim{{\xi}}    
\def\ra{$r$-adaptivity}
\def\dfem{$\partial$FEM}

\newcommand\anon[2]{{#2}} 

\title{Differentiable High-Order Mesh $r$-Adaptivity}
\author{
Ketan Mittal\thanks{Corresponding author, mittal3@llnl.gov.} \qquad
Julian Andrej \qquad
John Camier \qquad
Veselin Dobrev \\
Tzanio Kolev \qquad
Leonardo Molinari \qquad
Vladimir Tomov
}
\date{}

\makeatletter
\g@addto@macro\@thanks{\protect\footnotetext[0]{Lawrence Livermore National Laboratory, Livermore, CA, U.S.A.}}
\makeatother

\begin{document}
\newcommand\relatedversion{}

\maketitle

\begin{abstract}
\small\baselineskip=9pt
High-order mesh $r$-adaptivity typically entails minimizing a nonlinear functional using its gradients and, for Newton-based methods, its Hessian. These derivatives can be expensive to derive, must be updated with any changes to the functional, and errors in the derivatives directly impact the convergence of \ra. Additionally, derivatives of the \ra\ operator are required for consistent sensitivity analysis for inverse problems. This work develops a framework for differentiable \ra\ based on the Target-Matrix Optimization Paradigm \cite{TMOP2019SISC}.
Given a pointwise functional defining the local mesh-adaptation objective,
the proposed framework uses \dfem\ \cite{andrej2026scalable} to
automatically generate pointwise first- and second-order derivatives.
These local derivatives are then assembled into a global gradient vector
and a Hessian operator through the finite element operator decomposition.
The derivative kernels support portable CPU and GPU execution,
with partial assembly and tensor sum factorization enabling matrix-free application of the Hessian operator. The framework
also provides consistent sensitivities of the \ra\ operator through
automatically generated mixed derivatives and implicit differentiation
of the mesh stationarity equations. Numerical experiments show GPU solve times
comparable to those of hand-derived TMOP implementations on a benchmark problem, improvement in $r$-adaptivity convergence due to accurate computation of the derivatives, and demonstrate the use of mesh sensitivities for material-interface identification in an inverse heat-conduction problem.
\end{abstract}

\section{Introduction}
\label{sec:introduction}

High-order meshes are essential in the high-order finite element method (FEM)
for accurately modeling the geometry, achieving optimal convergence rates
on domains with curved boundaries/interfaces, preserving symmetry, and
aligning the mesh with the key features of the flow in moving mesh simulations
\cite{Dobrev2012}.
The quality of mesh elements, which encompasses their shape, size, and orientation, directly impacts the accuracy and stability of the numerical solution. For example, poor element shapes can
lead to ill-conditioned algebraic systems \cite{mittal2019mesh}, while
insufficient resolution can limit solution accuracy \cite{TMOP2021EWC}.
Mesh $r$-adaptivity improves the quality of an existing mesh by
moving its nodes without changing the mesh topology or the number of
degrees of freedom. The node movement can also be guided by the numerical
solution to improve local resolution and alignment with solution features.
Existing approaches for high-order mesh $r$-adaptivity include the Target-Matrix Optimization Paradigm (TMOP) \cite{TMOP2019SISC,TMOP2020CAF,TMOP2021EWC,TMOP2023CAD} and metric tensor-based adaptation \cite{frey2005anisotropic,alauzet2016decade,aparicio2023combining,marcon2019variational}.

These approaches formulate $r$-adaptivity as the minimization of a nonlinear functional of the mesh node coordinates, amongst other parameters.
Efficient and robust optimization then requires accurate first derivatives for gradient-based solvers such as the Broyden–Fletcher–Goldfarb–Shanno (BFGS) method, and both first and second derivatives for Newton's method \cite{TMOP2019SISC,aparicio2024globalized}.
Deriving and implementing these derivatives manually is time consuming and error prone. Introducing any modifications to the functional requires corresponding changes to the derivatives, which makes prototyping new formulations slow. The implementation effort increases further when supporting full assembly, partial assembly, and GPU execution, each of which may require a specialized implementation.

These challenges are further exacerbated when sensitivity analysis of the
mesh optimization operator itself is required. In this case, one must account
for the dependence of the optimized mesh on the inputs to the \ra\ operator,
and combine the resulting sensitivities with those of the upstream and downstream operators (e.g., a PDE solve) through the chain rule.
For example, sensitivity analysis of the \ra\ operator is needed for
PDE-constrained mesh optimization \cite{pde2026},
topology optimization with conformal meshes \cite{schmidt2024level},
and for differentiable arbitrary Lagrangian-Eulerian (ALE) hydrodynamics when
\ra\ is used during the remeshing phase.

The current work uses the differentiable finite element framework
\dfem\ of Andrej et al. \cite{andrej2026scalable} to compute consistent
first- and second-order derivatives of TMOP-based high-order mesh
$r$-adaptivity functionals. \dfem\ supports automatic differentiation (AD) of pointwise functions using dual numbers \cite{baydin2018automatic} or, preferably, Enzyme \cite{moses2020enzyme,moses2021enzymegpu},
and combines these derivatives in a scalable manner using the finite element operator decomposition. \dfem\ is available open-source through MFEM \cite{mfem2024,CEED2021}, a high-order finite element library.

Automatic differentiation (AD) is a well-established technique that has been implemented using various techniques. Some popular approaches include AD based on operator overloading in TinyAD \cite{schmidt2022tinyad}, symbolic
differentiation in ACORNS \cite{desai2022acorns} and
SymX \cite{fernandez2026symx}, source code transformation in
Tapenade \cite{hascoet2013tapenade}, and compiler-level transformation in
Enzyme \cite{moses2020enzyme,moses2021enzymegpu}. The current work leverages the Enzyme-based \dfem\ framework as it provides a scalable and efficient way to compute derivatives of nonlinear finite element operators that support partial assembly with tensor sum factorization on GPU-accelerated hardware.

A survey of the literature shows that the use of automatic differentiation for
mesh generation and adaptivity has been historically limited,
but it has gained popularity in recent times.
TinyAD has been used to compute gradients and Hessians of local mesh energies on CPUs \cite{schmidt2022tinyad}, which has been extended for locality-aware derivative evaluation on GPUs \cite{mahmoud2026locality}.
Bunge et al. \cite{bunge2024polygon} use TinyAD for polygon mesh smoothing
by minimizing the trace of the discrete Laplacian, with an additional penalty
to preserve surface geometry and sharp features.
Similarly, AD has also been used for high-order mesh curving of linear meshes \cite{schneider2021isogeometric} where quadratic poly-spline isogeometric analysis (IGA) bases are constructed for each element such that the initial boundary is approximated closely. AD is used in this mesh curving process to compute derivatives of an element degeneration penalty term that ensures that element quality remains high. Thakur and Nadarajah \cite{thakur2025adjoint} use AD in goal-oriented
implicit shock tracking, where the objective combines an adjoint-weighted
residual error indicator with a particular mesh-distortion penalty.

The emphasis of the present work is on automatic differentiation of a family of TMOP-based high-order $r$-adaptivity objectives. To this end, the \dfem\ framework's support for differentiating nonlinear finite element residuals \cite{andrej2026scalable} has been extended to support scalar objectives such as those appearing in TMOP. The resulting formulation supports geometry- and solution-driven adaptivity through combinations of mesh-quality metrics \cite{knupp2022geometric,knupp2024seventeen}, spatially varying targets and node limiting \cite{TMOP2020CAF}, and surface-fitting terms \cite{TMOP2023CAD,TMOP2021IMR}.
The derivative operators retain the partial assembly and portable CPU and GPU execution provided by \dfem, allowing an efficient matrix-free approach for \ra.
The proposed framework is demonstrated to be effective at generating accurate derivatives of TMOP functionals, and for sensitivity analysis of the \ra\ operator.

The remainder of the paper is organized as follows. Section 2 gives an overview of the TMOP-based $r$-adaptivity formulation and \dfem\ framework.
Section 3 describes how \dfem\ is used to compute the derivatives of the TMOP-based $r$-adaptivity functional, as well as the sensitivity associated with the $r$-adaptivity operator.
Section 4 presents numerical examples that illustrate the effectiveness and efficiency of the proposed approach. Finally, Section 5 concludes the paper with a summary of the main contributions and potential directions for future work.

\section{Preliminaries}
\label{sec:prelim}
This section describes \dfem\ and the TMOP functional for mesh quality improvement, along with two modifications, namely node limiting and mesh fitting, that are relevant to our applications of interest.

\subsection{High-order mesh and function representation}
\label{sec:mesh-representation}

In the finite element method, the domain $\Omega \subset \R^d$, $d\in\{2,3\}$, is discretized as a union of curved mesh elements, $\Omega^e$, $e=1,\dots,N_E$, each of order $p$.
A discrete representation of these elements is obtained by selecting a set of scalar basis functions $\{w_i\}_{i=1}^{N_p}$ on the reference element $\bar{\Omega}^e$.
In the case of tensor-product elements (quadrilaterals in 2D and hexahedra
in 3D), $N_p=(p+1)^d$, and the basis spans the space of all polynomials of
degree at most $p$ in each variable, and $\bar{\Omega}^e=[0,1]^d$.
Typically, Lagrange interpolants on the Gauss-Lobatto nodes of the reference
element are chosen for the $p$th-order basis functions. For a conforming
mesh, $C^0$ continuity is enforced by matching the nodal coordinates on
shared element boundaries.

Let $\vX$ denote the global vector of nodal coordinates, i.e., the mesh
\emph{degrees of freedom} (DOFs). The coordinates for element $\Omega^e$
are obtained by restricting $\vX$ to that element and arranging the entries
in a $d\times N_p$ matrix $\vx_e$. The $i$th column, $\vx_{e,i}$, contains the
coordinates of local node $i$. These coordinates define the position
function through the map $\Phi_e:\bar{\Omega}^e\to\Omega^e$:
\begin{equation}
\label{eq_x}
\begin{aligned}
\vx(\bar{\vx};\vX)&=\Phi_e(\bar{\vx};\vX)
  =\sum_{i=1}^{N_p}\vx_{e,i}w_i(\bar{\vx}),
\end{aligned}
\end{equation}
where $\bar{\vx}\in\bar{\Omega}^e$.

Mesh $r$-adaptivity updates $\vX$ while keeping the element connectivity,
polynomial order $p$, reference basis functions, and quadrature rule fixed.

\subsection{TMOP for Mesh Quality Improvement via $r$-adaptivity}
\label{sec:r-adaptivity}
The Jacobian of the mapping $\Phi_e$ at reference point
$\bar{\vx}\in\bar{\Omega}^e$ is defined as
\begin{equation}
\label{eq_A}
\begin{aligned}
  \left[A_e(\bar{\vx};\vX)\right]_{ab}
    &=\frac{\partial x_a(\bar{\vx};\vX)}{\partial\bar{x}_b}
    =\sum_{i=1}^{N_p}(\vx_{e,i})_a
      \frac{\partial w_i(\bar{\vx})}{\partial\bar{x}_b},
\end{aligned}
\end{equation}
where $a,b=1,\dots,d$. This matrix describes the local deformation relative
to the reference element. In \eqref{eq_A}, $x_a$ is the $a$th component of
the position function and $(\vx_{e,i})_a$ is the $a$th coordinate of the
$i$th node in element $\Omega^e$. Both the position function \eqref{eq_x} and the Jacobian \eqref{eq_A} depend linearly on the nodal coordinates.
The Jacobian matrix $A$ can impact the accuracy and computational cost of the solution \cite{mittal2019mesh}, which lends to the central idea of TMOP of optimizing the mesh to control the local Jacobian $A_{d \times d}$ in the mesh.

\subsubsection*{Target matrix}
The target matrix $W(\vx)\in\R^{d\times d}$ describes the desired local
Jacobian transformation from the reference element $\bar{\Omega}^e$ to an
ideal target element $\Omega^{e_t}$ that specifies the desired
element shape, size, and orientation.
Target construction is typically guided by the fact that any Jacobian matrix can be written as a composition of four geometric components \cite{knupp2022geometric}, namely volume, rotation, skewness, and aspect ratio:
\begin{equation}
\label{eq_W}
\begin{aligned}
W(\vx) &= \underbrace{\zeta(\vx)}_{\text{[volume]}} \circ
    \underbrace{R(\vx)}_{\text{[rotation]}} \circ
    \underbrace{Q(\vx)}_{\text{[skewness]}} \circ
    \underbrace{D(\vx)}_{\text{[aspect ratio]}}.
\end{aligned}
\end{equation}
In practice, the user may specify $W$ as a combination of any of the four fundamental components that they are interested in optimizing the mesh for.
Target construction is not the focus of this paper, and the reader is referred to \cite{TMOP2020CAF,TMOP2021EWC,knupp2019target} for detailed discussion on this topic.
The target matrix, or its components, may be defined analytically or through discrete finite element functions. In the latter scenario, their value at a given reference point depends on the current mesh configuration through the mapping $\Phi_e$:
\begin{equation}
\label{eq:target-dependence}
W_e(\bar{\vx};\vX)=W\big(\vx(\bar{\vx};\vX)\big).
\end{equation}
When targets are defined through discrete functions on the initial mesh, they must be remapped to the current nodal positions after each \ra\ iteration. This is typically done via general high-order field evaluation utilities such as FindPoints \cite{mittal2025general}.

Hereafter, we suppress the arguments $(\bar{\vx};\vX)$ of elementwise
quantities and use the subscript $q,e$ for evaluation at quadrature point
$\bar{\vx}_q$ in element $e$, e.g.,
$A_{q,e}=A_e(\bar{\vx}_q;\vX)$. Dependence on $\vX$ remains implicit unless
a particular mesh iterate must be specified.

\subsubsection*{Mesh quality metric}
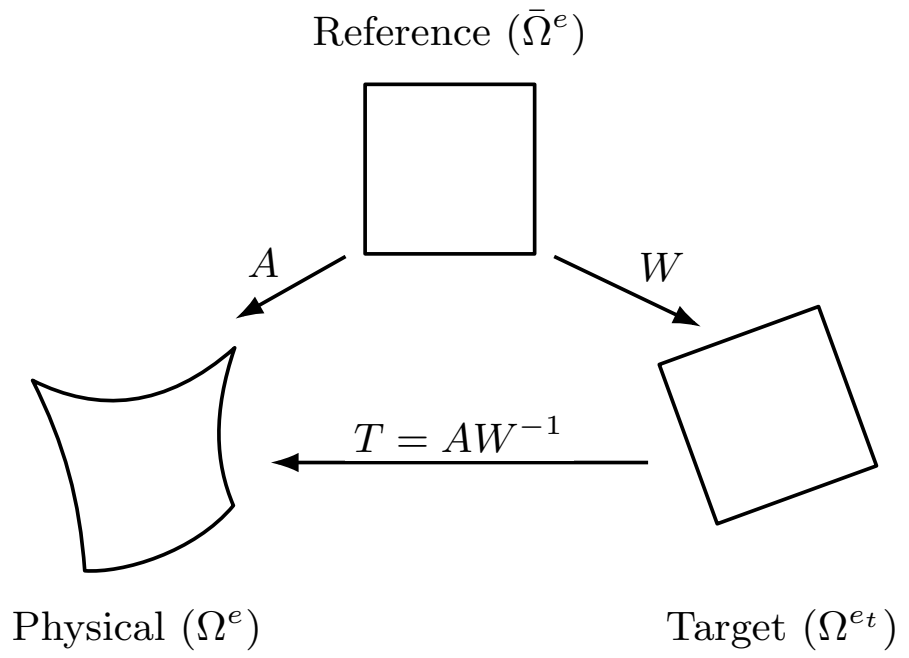
\begin{figure}[tb!]
\centerline{\resizebox{0.8\linewidth}{!}{\begingroup
\newcommand{\tmopReferenceLabel}{Reference ($\bar{\Omega}^e$)}
\newcommand{\tmopPhysicalLabel}{Physical ($\Omega^e$)}
\newcommand{\tmopTargetLabel}{Target ($\Omega^{e_t}$)}
\newcommand{\tmopPhysicalMapLabel}{$A$}
\newcommand{\tmopTargetMapLabel}{$W$}
\newcommand{\tmopRelativeMapLabel}{$T=AW^{-1}$}

\begin{tikzpicture}[
  x=1cm, y=1cm,
  font=\small,
  element/.style={draw=black, line width=0.8pt, line join=round},
  jacobian/.style={-{Latex[length=2.2mm,width=1.5mm]}, line width=0.8pt},
  map label/.style={fill=white, inner sep=2pt}
]
  \draw[element] (2.85,2.62) rectangle (4.15,3.92);
  \node[anchor=south] at (3.50,4.02) {\tmopReferenceLabel};

  \draw[element]
    (0.30,1.65)
      .. controls (0.80,1.40) and (1.32,1.43) .. (1.85,1.90)
      .. controls (1.68,1.40) and (1.70,1.02) .. (1.84,0.69)
      .. controls (1.54,0.33) and (0.99,0.17) .. (0.70,0.19)
      .. controls (0.66,0.73) and (0.55,1.15) .. (0.30,1.65)
    -- cycle;
  \node[anchor=north] at (1.10,0.02) {\tmopPhysicalLabel};

  \begin{scope}[shift={(5.55,0.55)}, rotate=20]
    \draw[element] (0,0) rectangle (1.30,1.30);
  \end{scope}
  \node[anchor=north] at (6.05,0.02) {\tmopTargetLabel};

  \draw[jacobian] (2.70,2.60) --
    node[map label, above left] {\tmopPhysicalMapLabel} (1.85,2.12);
  \draw[jacobian] (4.30,2.60) --
    node[map label, above right] {\tmopTargetMapLabel} (5.43,2.06);
  \draw[jacobian] (5.02,1.02) --
    node[map label, above] {\tmopRelativeMapLabel} (2.12,1.02);
\end{tikzpicture}
\endgroup}}
\vspace{-5mm}
\caption{Schematic of the TMOP Jacobian matrices for the reference element
$\bar{\Omega}^e$, physical element $\Omega^e$, and target element
$\Omega^{e_t}$.}
\label{fig:tmop}
\vspace{-3mm}
\end{figure}

With the target matrix defined, a mesh quality metric is chosen to measure the deviation between the current and target Jacobian.
The mesh quality metric $\mu(T)$, evaluated at $T=AW^{-1}$ (Figure \ref{fig:tmop}), compares $A_e$ and $W_e$ in terms of one or more of their
geometric parameters.
For example, in two dimensions,
$\mu_{2,s}(T)=\frac{\lVert T\rVert_F^2}{2\tau}-1$ is a \emph{shape} metric\footnote{The metric subscript follows the numbering in \cite{Knupp2020, knupp2022geometric}.}
that depends on the skewness and aspect ratio components, but
is invariant to orientation/rotation and volume. Here, $\lVert T\rVert_F$ denotes the Frobenius norm and $\tau=\det T$.
Similarly, $\mu_{36,vos}(T)=\frac{\lVert T-I\rVert_F^2}{\tau}$ is a
\emph{shape}$+$\emph{size}$+$\emph{orientation} metric
that penalizes differences in each of the four geometric parameters.
For $\tau>0$, its minimum value is zero and is attained only when
$T=I$, i.e., when $A=W$, where $I$ denotes the identity matrix.
These different mesh quality metrics allow the user to set up TMOP-based $r$-adaptivity according to their simulation needs.

\subsubsection*{TMOP objective function}

With the target matrix and mesh quality metric defined, $r$-adaptivity is done by minimizing the global mesh-quality objective:
\begin{equation}
\label{eq:F_full}
F_\mu(\vX)=\sum_{e=1}^{N_E}\int_{\bar{\Omega}^e}
  \omega_e\,\mu(T_e)\,\dd\bar{\vx},
\end{equation}
where $\omega_e=\det(W_e)>0$ accounts for the target size associated with element $\Omega^{e_t}$. The integral \eqref{eq:F_full} is
evaluated using a standard quadrature-based approach:
\begin{equation}
\label{eq:F_full_quadrature}
F_\mu(\vX)\approx\sum_{e=1}^{N_E}\sum_{q=1}^{N_q^e}
\beta_q\,\omega_{q,e}\,\mu(T_{q,e}),
\end{equation}
where $N_q^e$ is the number of quadrature points in element $e$ and $\beta_q$
is the fixed reference quadrature weight, with
$T_{q,e}=A_{q,e}W_{q,e}^{-1}$ and $\omega_{q,e}=\det(W_{q,e})$.
Minimizing \eqref{eq:F_full_quadrature} results in node movement such that
the local Jacobian transformation $A$ resembles the target transformation $W$ as closely as possible at each point, in terms of the geometric parameters enforced by the mesh quality metric.

\subsubsection*{Node limiting}
\label{sec:node-limiting}

In applications such as arbitrary Lagrangian-Eulerian (ALE) simulations, the TMOP objective is modified to limit node movement by adding a displacement penalty \cite{TMOP2020CAF},
\begin{equation}
\label{eq:node-limiting}
F_{\mathrm{lim}}(\vX)=\sum_{e=1}^{N_E}
\int_{\bar{\Omega}^e}\frac{\|\vx-\vx_0\|^2}{2c_{\mathrm{lim}}^2}\,\omega_e\,\dd\bar{\vx}.
\end{equation}
Here, $\vx=\vx(\bar{\vx};\vX)$, $\vx_0(\bar{\vx})$ is the fixed initial position function, $\omega_e$ is the size associated with the target configuration, and $c_{\mathrm{lim}}>0$ is a penalty parameter based on the local limiting distance.
Minimizing $F_\mu+F_{\mathrm{lim}}$ balances mesh quality improvement with keeping nodes in the vicinity of their original positions. This penalty is especially useful in ALE simulations when the mesh deforms with the material flow and improving mesh quality while simultaneously limiting node movement is necessary to avoid numerical diffusion during the subsequent remap phase. We consider $c_{\mathrm{lim}}$ to be constant in \eqref{eq:node-limiting}, but it can vary spatially if the application requires.
In practice, $c_{\mathrm{lim}}$ is automatically chosen based on node displacement observed during the Lagrangian phase prior to the mesh improvement step; see \cite{TMOP2020CAF} for further details.

\subsubsection*{Mesh fitting to implicit geometries}

The TMOP-based framework also enables mesh generation via adaptation for implicit geometries where the domain boundaries/interfaces are represented by the zero isocontour of a level-set function. In this case, a subset of mesh nodes is weakly enforced to lie on the zero isosurface of the level-set function by adding a penalty term to the objective,
\begin{equation}
\label{eq:surface-fitting}
F_{\mathrm{surf}}(\vX)=c_{\mathrm{surf}} \sum_{s \in \mathcal{S}}
\sigma\big(\vx_s(\vX)\big)^2,
\end{equation}
where $c_{\mathrm{surf}}\geq0$ is the penalization weight, $\mathcal{S}$ is the set of nodes to be aligned to the level set (e.g., mesh nodes discretizing the domain boundary or multimaterial interface), and $\sigma(\vx_s(\vX))$ is the level-set function evaluated at the coordinate vector $\vx_s(\vX)$ of node $s$. Minimizing $F_\mu+F_{\mathrm{surf}}$ weakly enforces the selected nodes to lie on the implicit surface while still improving mesh quality, and produces a boundary- or interface-fitted mesh. A detailed discussion of this formulation along with ways to make it robust for practical computations is in \cite{TMOP2021IMR,TMOP2023CAD,mittal2023high}.

\subsection{TMOP Solver}
\label{sec:tmop-solver}
The global \ra\ objective,
$F(\vX)=F_\mu(\vX)+F_{\mathrm{lim}}(\vX)+F_{\mathrm{surf}}(\vX)$,
is minimized iteratively using Newton's method:
\begin{equation}
\label{eq:r_adaptivity_solve}
\vX_{k+1} = \vX_k - \alpha_k \mathcal{H}^{-1}(\vX_k)\mathcal{J}(\vX_k),
\end{equation}
where $\vX_k$ is the global vector of nodal coordinates at the $k$th Newton
iteration, $\alpha_k$ is the step length determined by a line-search
procedure, and $\mathcal{J}(\vX)=\nabla_{\vX}F(\vX)$ and
$\mathcal{H}(\vX)=\nabla_{\vX}^2F(\vX)$ are the gradient and Hessian of the
total objective, respectively.
The line-search procedure is described in prior work \cite{TMOP2019SISC,TMOP2020CAF,TMOP2021EWC}.
Newton iterations \eqref{eq:r_adaptivity_solve} are done until the line-search criteria cannot be satisfied or the gradient norm relative to its initial value satisfies $\|\mathcal{J}(\vX_k)\|_2/\|\mathcal{J}(\vX_0)\|_2\leq\varepsilon$,
where $\vX_0$ denotes the initial nodal coordinates and $\varepsilon=10^{-10}$ by default.

In the current work, the focus is on automating computation of the TMOP derivatives, i.e., the gradient $\nabla_{\vX}F(\vX)$ and the Hessian  $\nabla_{\vX}^2F(\vX)$ using \dfem.

\subsection{\dfem: Differentiable Finite Element Method}
\label{sec:dfem}

In this work, differentiation is enabled through the \dfem\ interface in MFEM, previously described in~\cite{mfem2024,andrej2026scalable}.
We briefly introduce its general approach to differentiability and then discuss some features that support the present work. 

\begin{figure*}[!t]
\centerline{\includegraphics[width=0.8\textwidth]{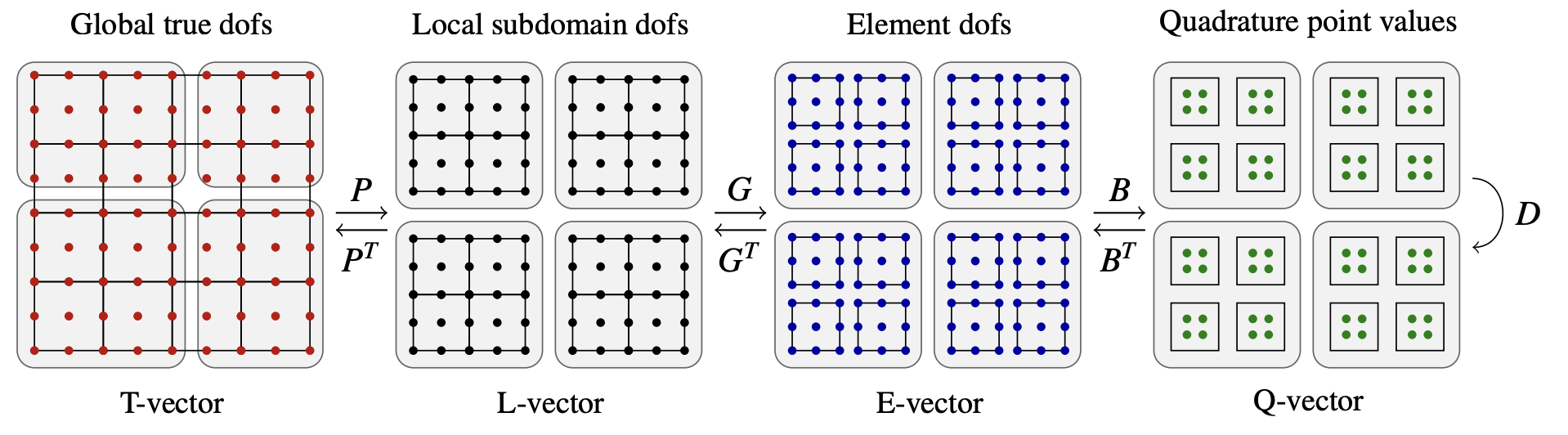}}
\vspace{-5mm}
\caption{Finite element operator decomposition. The operators $P$, $G$,
and $B$ map global true degrees of freedom to local subdomain degrees of
freedom, element degrees of freedom, and quadrature-point values,
respectively. The pointwise operator $D$ acts at quadrature points, and
the transpose operators accumulate the contributions back to the global
degrees of freedom.}
\label{fig:operator-decomposition}
\end{figure*}

Building on the existing MFEM infrastructure, \dfem\ leverages the finite element operator decomposition \cite{mfem2024}, illustrated in Figure \ref{fig:operator-decomposition} and written here for a linear operator $\mathcal{A}$ as:
\begin{equation}
  \mathcal{A} = P^{\mathsf T}G^{\mathsf T}B^{\mathsf T}D BGP,
\end{equation}
where $P$ maps global true degrees of freedom to local subdomain degrees of freedom, $G$ restricts them to individual elements, and $B$ evaluates fields and their reference derivatives at quadrature points. The operator $D$ contains the pointwise physics, geometric transformations, and quadrature weights.
For fixed mesh connectivity, reference finite element spaces, and quadrature rules, $P$, $G$, and $B$ are \emph{purely topological} and \emph{linear} maps, which do not depend on the physics or the solution.

This allows AD to be restricted to the pointwise kernel, without differentiating additional layers such as MPI communication, time integration or nonlinear solvers.
Localizing AD in this way further limits its computational and memory overhead and minimizes changes to existing code.

Following the same decomposition, a nonlinear residual operator can be written as
\begin{equation}
  \mathcal{R}(u)=P^{\mathsf T}G^{\mathsf T}B^{\mathsf T}D(u_q),
  \qquad u_q=BGPu.
\end{equation}
Its Jacobian follows directly from the chain rule:
\begin{equation}
  J_{\mathcal{R}}(u)=P^{\mathsf T}G^{\mathsf T}B^{\mathsf T}
  \frac{\partial D}{\partial u_q}(u_q)BGP.
\end{equation}
Since $D$ acts pointwise, $\frac{\partial D}{\partial u_q}$ is block diagonal, with one small block per quadrature point. Each block can be evaluated locally and independently.
In practice, the user supplies the pointwise kernel, while \dfem\ handles differentiation and composes the resulting derivatives with MFEM's existing basis evaluations and gather-scatter operations.
For a scalar functional, the gradient is a nonlinear residual
operator, and its Jacobian is the Hessian.

Thanks to this architecture, the \dfem\ interface can run efficiently on both CPUs and GPUs, enabling portability across CUDA and HIP backends available in MFEM.

Quadrature-point kernels are differentiated using Enzyme~\cite{moses2020enzyme,moses2021enzymegpu}, an LLVM compiler plugin that operates on intermediate representation (IR)\footnote{A
native MFEM dual-number implementation is also available as an alternative AD backend.}.
This allows Enzyme to support multiple source languages targeting LLVM and to benefit from compiler optimizations when generating derivative code.

Enzyme supports both forward- and reverse-mode AD. For residual Jacobians, we favor forward mode because of its reduced memory overhead compared to reverse mode, particularly on GPUs.

\subsubsection*{Assembly levels}
Analogously to native MFEM operators, \dfem\ allows different assembly levels for derivatives: matrix-free evaluation (MF), partial assembly (PA), and full matrix assembly (FA).
For tensor-product elements, the basis evaluations retain MFEM's sum-factorization structure.

In the MF case, the primal and derivative actions are evaluated on the fly; PA instead separates the initial setup stage, evaluating and caching quadrature-point derivative data at a given state, from the subsequent application.
This is useful when a linear solver repeatedly applies the derivative operator at the same linearization state, as in the present work.
Finally, Jacobians can still be assembled into a global
sparse matrix for FA if needed.

\subsubsection*{Registered derivatives in Enzyme}
A user may wish to provide a custom derivative for a specific part of the computation (e.g., when the derivative cannot be computed using pointwise data).
Enzyme supports this through derivative registration, which supplies a differentiation rule in place of differentiating the underlying implementation.
The registered rule is then incorporated into the overall derivative through the chain rule.
This mechanism is used for discrete target and level-set fields, as
described in Sections \ref{sec:mesh-quality}
and \ref{sec:fitting-qfunction}.

\section{Methodology}
\label{sec:methodology}

We extend the residual-based differentiation framework of
\dfem\ to scalar TMOP objectives, obtaining gradients and Hessians from
user-provided q-functions, i.e., local functions evaluated at quadrature
points. These local derivatives are assembled into global operators using
the finite element operator decomposition
(Figure \ref{fig:operator-decomposition}). We first formulate the q-functions
for the mesh-quality, node-limiting, and mesh-fitting objectives introduced
in Section \ref{sec:prelim} and describe their differentiation, including
the chain rule and custom derivative rules for discrete target and
level-set fields. We then show how the resulting derivative operators
provide consistent sensitivities of the converged mesh with respect to
the inputs to the $r$-adaptivity operator.

\subsection{Mesh-quality q-function}
\label{sec:mesh-quality}

For \dfem, the mesh-quality objective \eqref{eq:F_full_quadrature} is
expressed as a q-function that takes the local position $\vx_{q,e}$ \eqref{eq_x} and reference Jacobian $A_{q,e}$ \eqref{eq_A} as inputs:
\begin{equation}
\label{eq:quality-qfunction}
f_{\mu,q,e}(\vx,A)
 =\beta_q\,\omega(\vx)\,\mu\big(AW(\vx)^{-1}\big),
\end{equation}
where $\omega(\vx)=\det W(\vx)$. The position is provided as an input because it may be needed for defining spatially varying target matrices $W(\vx)$, and the Jacobian $A_{q,e}$ is included explicitly as an input because it cannot be inferred solely from the local position $\vx_{q,e}$.

\subsubsection*{Gradient computation}
\label{sec:gradient-computation}

Suppressing the subscripts $q,e$, Enzyme automatically generates the partial derivatives of $f_\mu(\vx,A)$ with respect to each \emph{active} input while holding the other inputs fixed,
\begin{equation}
\label{eq:qfunction-partials}
\left.\frac{\partial f_\mu}{\partial\vx}\right|_A\in\R^d,
\qquad
\left.\frac{\partial f_\mu}{\partial A}\right|_{\vx}\in\R^{d\times d}.
\end{equation}

Since the dependence on the position $\vx$ enters
entirely through the target matrix $W(\vx)$, the chain rule gives
\begin{equation}
\label{eq:target-position-chain-rule}
\left.\frac{\partial f_\mu}{\partial\vx}\right|_A
 =\frac{\partial f_\mu}{\partial W}:
   \frac{\partial W}{\partial \vx},
\end{equation}
where $W$ enters $f_\mu$ through both
$\det W$ and $AW^{-1}$.

For an analytic target, $W(\vx)$ is constructed within the q-function from analytic expressions supplied by the user. Enzyme differentiates
the complete composition, including the target construction,
determinant, inverse, and metric evaluation. No explicit target
derivatives are required from the user.

For targets defined by discrete fields, one or more geometric
components of $W$ are prescribed using grid functions. For example,
the target size can be chosen based on the magnitude of the solution
gradient \cite{TMOP2020CAF}. In this discrete target case, the related grid functions are also provided as \emph{active} inputs to the q-function.
However, these local values provided to the q-function do not expose their spatial dependence and the grid function gradients must be computed outside using standard finite element techniques and provided to Enzyme through custom differentiation rules. Enzyme then handles the chain rule for associated derivative computation.

Using the finite element operator decomposition from
Section \ref{sec:dfem}, the quadrature-point inputs and global gradient
are computed with partial assembly \cite{CEED2021,camier2023accelerating}:
\begin{equation}
\label{eq:pa-quadrature-inputs}
\vx_Q=B_0GP\vX,
\qquad
\boldsymbol A_Q=B_1GP\vX.
\end{equation}
\begin{equation}
\label{eq:pa-nodal-gradient}
\begin{aligned}
\mathcal J_\mu(\vX)&=\frac{\partial F_\mu}{\partial\vX}\\
 &=P^{\mathsf T}G^{\mathsf T}\left(B_0^{\mathsf T} \left\{\frac{\partial f_{\mu,q,e}}{\partial\vx}\right\}
 +B_1^{\mathsf T} \left\{\frac{\partial f_{\mu,q,e}}{\partial A}\right\}\right).
\end{aligned}
\end{equation}
Here, $B_0$ and $B_1$ are the basis value and reference gradient
operators, respectively. In \eqref{eq:pa-nodal-gradient}, the first term accounts for dependence
through the quadrature-point positions and the spatially varying
target, while the second term accounts for dependence through the
reference Jacobians.

\subsubsection*{Hessian evaluation}
\label{sec:hessian-evaluation}

Second derivatives are computed with Enzyme using a \emph{forward-over-reverse} approach. Reverse-mode AD first generates the local q-function gradient; forward mode then differentiates this gradient with respect to its active inputs to compute Hessian data or actions.\footnote{The dual-number backend uses nested dual numbers to produce second derivatives.} At each quadrature point, this gives the blocks
\begin{equation}
\label{eq:qfunction-hessian}
D_{q,e}=
\begin{bmatrix}
\displaystyle\frac{\partial^2 f_\mu}{\partial\vx^2}
&
\displaystyle\frac{\partial}{\partial A}
 \left(\frac{\partial f_\mu}{\partial\vx}\right)
\\[2mm]
\displaystyle\frac{\partial}{\partial\vx}
 \left(\frac{\partial f_\mu}{\partial A}\right)
&
\displaystyle\frac{\partial^2 f_\mu}{\partial A^2}
\end{bmatrix}_{q,e}.
\end{equation}
This procedure generates both the pure and mixed derivative blocks in \eqref{eq:qfunction-hessian}.

Let $B=[B_0^{\mathsf T}\;B_1^{\mathsf T}]^{\mathsf T}$, and let $D(\vX)$ apply the blocks $D_{q,e}$ independently at each quadrature point. The action of the mesh-quality Hessian $\mathcal H_\mu=\partial^2F_\mu/\partial\vX^2$ on a nodal vector $\boldsymbol v$ is
\begin{equation}
\label{eq:pa-hessian-action}
\mathcal H_\mu(\vX)\boldsymbol v
 =P^{\mathsf T}G^{\mathsf T}B^{\mathsf T}D(\vX)BGP\boldsymbol v.
\end{equation}
Combining the derivatives of the mesh-quality, node-limiting, and
mesh-fitting terms gives the total gradient $\mathcal J$ and Hessian
$\mathcal H$ defined in Section \ref{sec:tmop-solver}.
The total gradient $\mathcal J(\vX)$ is a nonlinear operator evaluated at the supplied mesh state. Its Jacobian, the Hessian $\mathcal H(\vX_k)$, is constructed at a fixed Newton linearization state and acts linearly on perturbations. Hessians support the assembly levels described in Section \ref{sec:dfem}; we use partial assembly. At each Newton linearization, the local Hessian blocks are evaluated once and cached for reuse throughout the Krylov solve for computing $\Delta\vX_k$ in $\mathcal H(\vX_k)\Delta\vX_k=\mathcal J(\vX_k)$
\eqref{eq:r_adaptivity_solve}.

The same differentiation procedure also yields mixed derivatives with respect to different fields or parameters supplied as active inputs. Note that for analytic targets, Enzyme differentiates the target construction and energy evaluation a second time. For targets constructed from grid functions, the custom differentiation rules also supply the fields' spatial Hessians, computed externally using finite element techniques. Enzyme then handles the chain rule for associated Hessian computation.

\subsection{Node-limiting q-function}
\label{sec:node-limiting-qfunction}

The q-function for the node-limiting objective \eqref{eq:node-limiting} is
\begin{equation}
\label{eq:limiting-qfunction}
f_{\mathrm{lim},q,e}(\vx)
 =\frac{\beta_q}{2c_{\mathrm{lim}}^2}\,\omega(\vx)
   \|\vx-\vx_{0,q,e}\|^2,
\end{equation}
where $\vx_{0,q,e}$ is the initial position evaluated at the same
reference quadrature point. With $\vx_{0,q,e}$ and
$c_{\mathrm{lim}}$ fixed during differentiation, $\vx$ is the only
active input. Enzyme therefore computes the local gradient
$\partial f_{\mathrm{lim}}/\partial\vx$ and Hessian
$\partial^2 f_{\mathrm{lim}}/\partial\vx^2$ directly.
This includes the position dependence of $\omega(\vx)$, using the
same analytic target expressions or custom field rules as for mesh
quality in Section \ref{sec:mesh-quality}.

Since the q-function has no dependence on $A$, the global gradient
uses only the $B_0^{\mathsf T}$ contraction in \eqref{eq:pa-nodal-gradient}.
The Hessian action follows \eqref{eq:pa-hessian-action} with $B=B_0$
and the local Hessians $\partial^2 f_{\mathrm{lim}}/\partial\vx^2$
as quadrature blocks.

\subsection{Mesh-fitting q-function}
\label{sec:fitting-qfunction}

For mesh fitting to implicit geometries, we use an integration rule $\mathcal N_e$ whose
points coincide with the element's mesh nodes. Let $s(q,e)$ denote
the global node corresponding to quadrature point $q$ and $n_s$ the
number of elements sharing node $s$ across all MPI ranks. The
q-function is
\begin{equation}
\label{eq:fitting-qfunction}
\begin{aligned}
f_{\mathrm{surf},q,e}(\vx)&=\widehat c_{q,e}\,\sigma(\vx)^2,\\
\widehat c_{q,e}&=
\begin{cases}
c_{\mathrm{surf}}/n_{s(q,e)},&s(q,e)\in\mathcal S,\\
0,&s(q,e)\notin\mathcal S.
\end{cases}
\end{aligned}
\end{equation}
Thus unmarked nodes contribute zero, and the factor $1/n_s$ prevents
multiple counting of shared nodes. The nodal quadrature sum
\begin{equation}
\label{eq:fitting-quadrature-sum}
F_{\mathrm{surf}}(\vX)
 =\sum_e\sum_{q\in\mathcal N_e}
 f_{\mathrm{surf},q,e}(\vx_{q,e})
\end{equation}
recovers \eqref{eq:surface-fitting}, with no additional quadrature
weight or Jacobian determinant.

With $\widehat c_{q,e}$ fixed, differentiating the fitting q-function
with respect to $\vx$ requires the spatial derivatives of $\sigma$.
For a discrete level set, these derivatives cannot be inferred from
the sampled value of $\sigma$ alone. As for discrete targets in
Section \ref{sec:mesh-quality}, we evaluate the spatial gradient
$\nabla_{\vx}\sigma$ and Hessian $\nabla_{\vx}^2\sigma$ using the
existing TMOP routines and provide them to Enzyme through custom
differentiation rules. Enzyme then applies the chain rule to compute
the local gradient $\partial f_{\mathrm{surf}}/\partial\vx$ and
Hessian $\partial^2 f_{\mathrm{surf}}/\partial\vx^2$.

\subsection{Sensitivity of the $r$-adaptivity operator}
\label{sec:radaptivity-sensitivity}

TMOP-based $r$-adaptivity is used in PDE-constrained mesh
optimization \cite{pde2026}, conformal meshing for level-set shape
optimization \cite{schmidt2024level}, and mesh relaxation in arbitrary
Lagrangian-Eulerian (ALE) hydrodynamics \cite{tomov2018}.
Consistent sensitivities of these simulations must account for the
dependence of the adapted mesh on the input parameters.

Let $\vX$ denote the free mesh coordinates and $\vtheta$ the input parameters defining targets or fitting interfaces. Parameter-independent constrained nodes (e.g., fixed boundary nodes) have zero sensitivity.
The converged stationary mesh $\vX^*(\vtheta)$ for the total objective $F(\vX,\vtheta)$ satisfies
\begin{equation}
\label{eq:mesh-stationarity}
\mathcal J\big(\vX^*(\vtheta),\vtheta\big)=\boldsymbol 0,
\qquad
\mathcal J=\frac{\partial F}{\partial\vX}.
\end{equation}

Differentiating \eqref{eq:mesh-stationarity} with respect to
$\vtheta$ gives
\begin{align}
\frac{\partial \mathcal J}{\partial \vX} \frac{\partial \vX^*}{\partial \vtheta} + \frac{\partial \mathcal J}{\partial \vtheta} &= 0, \\
\mathcal H_* \frac{\partial \vX^*}{\partial \vtheta} &= -\frac{\partial \mathcal J}{\partial \vtheta},
\end{align}
where $\mathcal H_*$ is the Hessian $\frac{\partial \mathcal J}{\partial \vX}$ evaluated at the converged mesh.
Assuming that $F$ is twice continuously differentiable near
$(\vX^*,\vtheta)$ and $\mathcal H_*$ is nonsingular, the sensitivity is
\begin{equation}
\label{eq:r-adaptivity-sensitivity}
\frac{\partial \vX^*}{\partial \vtheta} = -\mathcal H_*^{-1}\frac{\partial \mathcal J}{\partial \vtheta}.
\end{equation}
All partial derivatives are evaluated at $(\vX^*,\vtheta)$, treating
$\vX$ and $\vtheta$ as independent arguments of $F$.
With $\vtheta$ provided to the q-function as an \emph{active} input,
the mixed derivative $\partial\mathcal J/\partial\vtheta$ is also computed
using \dfem.

\subsection{Code Maintainability}
\label{sec:maintainability}

The q-function representation simplifies prototyping of mesh-quality
metrics and variations of the adaptation objective.
For example, a 2D ($\mu_{2,s}$) and 3D ($\mu_{303,s}$) shape-metric q-function is:
\begin{lstlisting}[mathescape=true,basicstyle=\small\ttfamily,
                   columns=fullflexible,keepspaces=true]
const auto W = target(x);
const auto T = A * inv(W);
const auto mu = $\frac{\lVert T\rVert_F^2}{d\,(\det T)^{2/d}}$ - 1.0;
f = w_q * det(W) * mu;
\end{lstlisting}
Here, \texttt{target(x)} is a separate user-specified routine for evaluating the target matrix that is common to all metrics $\mu(T)$,
\texttt{w\_q} is the reference quadrature weight provided by \dfem,
and $\lVert T\rVert_F$, $\det T$, and $inv()$ are \emph{factory} functions for computing the Frobenius norm, the determinant, and the matrix inverse, respectively.
The same q-function implementation can be conveniently modified for a different metric, e.g., for $\mu_{77}$ by replacing the metric definition with
\begin{lstlisting}[mathescape=true,basicstyle=\small\ttfamily,
                   columns=fullflexible,keepspaces=true]
const auto mu = $0.5 \cdot (\det T - \frac{1}{\det T})^2$;
\end{lstlisting}

The objective value, gradient, Hessian, and requested mixed derivatives
are generated using this implementation so it significantly reduces the
number of lines of code that the user must write and maintain for derivatives.
The same q-function implementation also uses the CPU and GPU execution
paths and retains partial assembly and tensor sum factorization provided by \dfem\ and MFEM.

\section{Numerical Results}
\label{sec:results}

This section presents various numerical experiments that demonstrate the effectiveness of the \dfem-based approach for portable CPU and GPU execution,
for improving convergence of \ra\ by accurately computing derivatives of the nonlinear TMOP functional, and for sensitivity analysis of the \ra\ operator for material-interface identification in an inverse heat-conduction problem.
All \dfem-based mesh-optimization runs in this work use partial assembly for the Hessian operator, on both CPUs and GPUs.

\subsection{Kershaw benchmark on CPU and GPU}
\label{sec:kershaw-results}

This benchmark assesses the performance impact of using the \dfem\ framework
by comparing mesh-optimization solve times with those of classic TMOP using
hand-derived derivatives on CPUs and GPUs. We consider the three-dimensional
Kershaw benchmark \cite{camier2023accelerating} for element shape control with linear, quadratic, and cubic meshes.

\begin{figure}[tbp]
\begin{center}
$\begin{array}{cc}
    \includegraphics[width=0.4\linewidth]{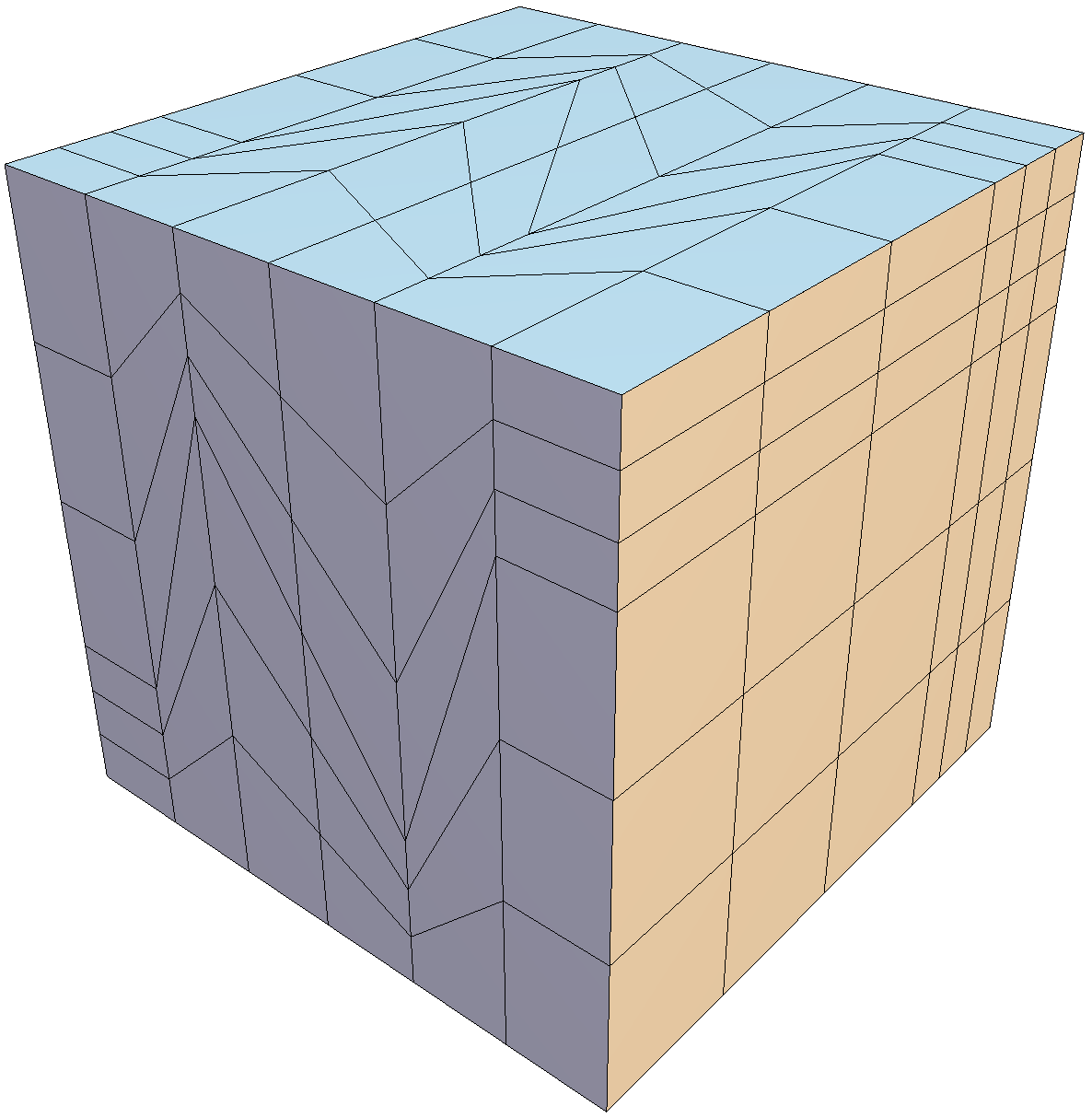} &
    \includegraphics[width=0.4\linewidth]{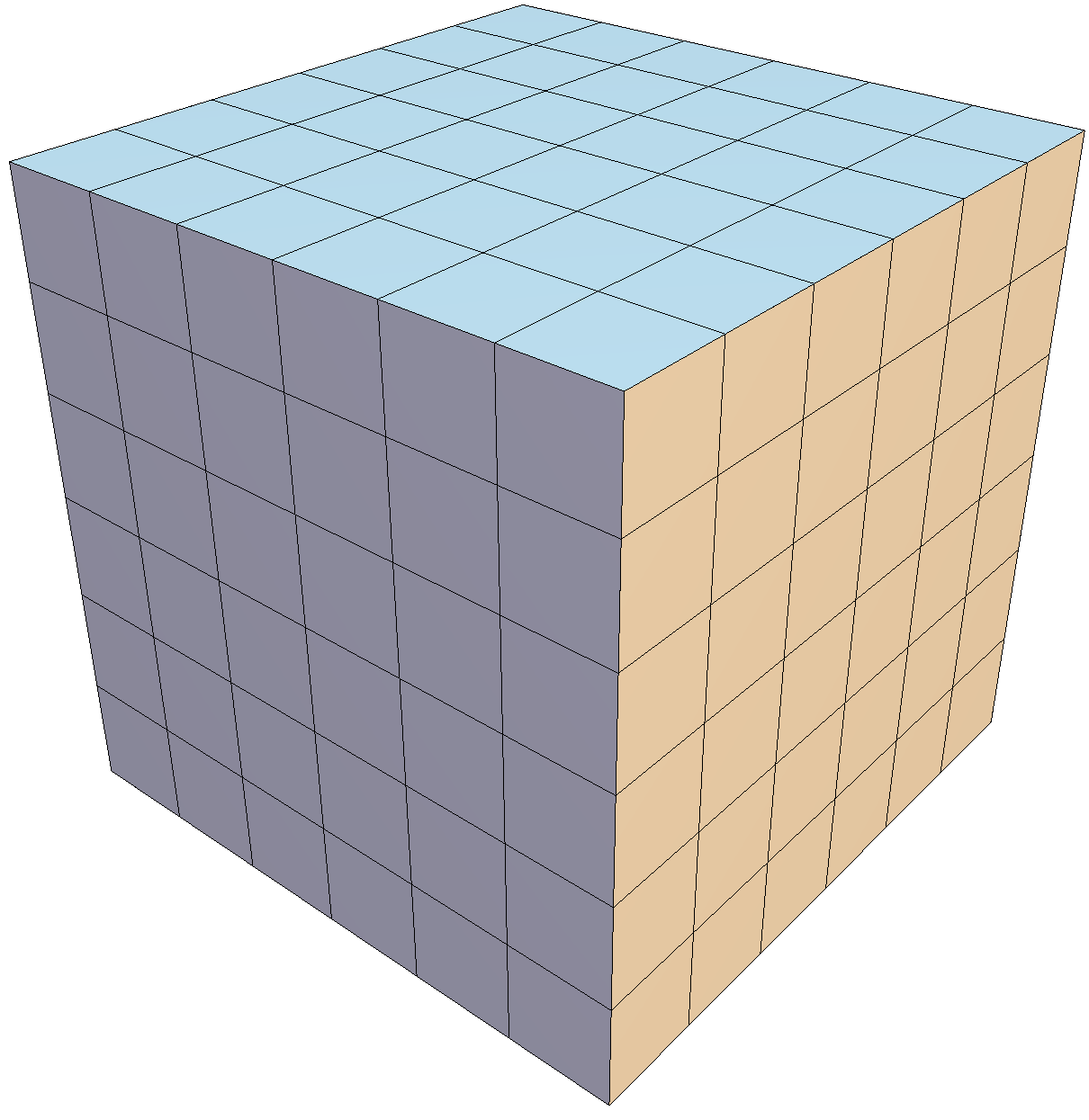} \\
    \textrm{(a)}\,\,\varepsilon_x=\varepsilon_z=0.5 &
    \textrm{(b)}\,\,\varepsilon_x=\varepsilon_z=1
\end{array}$
\end{center}
\vspace{-5mm}
\caption{Kershaw benchmark: (a) initial deformed mesh and (b) uniform
reference configuration recovered by optimization. The
$6\times6\times6$ meshes are shown for clarity; timing experiments
use three uniform refinements, giving $110{,}592$ elements.}
\label{fig:kershaw-meshes}
\end{figure}

The initial mesh is obtained by deforming a uniform Cartesian mesh
using the Kershaw transformation \cite{kershaw1981differencing, camier2023accelerating}, with anisotropy parameters
$\varepsilon_x=\varepsilon_z=0.5$.
Figure \ref{fig:kershaw-meshes} shows the initial configuration and
the uniform reference configuration, corresponding to
$\varepsilon_x=\varepsilon_z=1$, on a $6\times6\times6$ mesh.
For the timing experiments, a mesh with $48^3=110,592$ hexahedra is used.

In each case, the mesh with Kershaw deformation is optimized using the shape metric, $\mu_{303,s}(T)=\frac{\lVert T\rVert_F^2}{3\,\tau^{2/3}}-1$
to recover a uniform Cartesian mesh.
All the experiments are done on Matrix, a Livermore Computing system at Lawrence Livermore National Laboratory \cite{lc-matrix}. Each hardware node of Matrix contains two 56-core Intel Xeon Platinum 8480+ processors and four NVIDIA H100 GPUs.
The CPU runs use 56 CPU cores, while the GPU runs use one CPU core and one GPU on the same node. Table \ref{tab:kershaw-timings} reports the solve time $t$ and the number of Newton iterations $N$.

\begin{table}[t]
\centering
\caption{Kershaw benchmark timings on $110,592$ elements on Matrix comparing CPU and CPU+GPU performance for the \dfem-based framework in comparison to classic TMOP.}
\label{tab:kershaw-timings}
\vspace{2pt}
\setlength{\tabcolsep}{5pt}
\renewcommand{\arraystretch}{1.15}
\begin{tabular}{|c|c|c|c|c|}
\hline
& \multicolumn{2}{c|}{Solve time $t$ (s)}
  & \multicolumn{2}{c|}{Newton iterations $N$} \\
\hline
Order
  & Classic & \dfem\ & Classic & \dfem\ \\
\hline
\multicolumn{5}{|c|}{\textbf{56 CPU cores}} \\
\hline
1 & 38.97 & 58.66 & 10 & 10 \\
\hline
2 & 81.45 & 144.03 & 13 & 13 \\
\hline
3 & 257.57 & 529.89 & 21 & 22 \\
\hline
\multicolumn{5}{|c|}{\textbf{1 CPU core + 1 GPU}} \\
\hline
1 & 7.73 & 9.26 & 10 & 10 \\
\hline
2 & 15.93 & 15.49 & 13 & 13 \\
\hline
3 & 46.21 & 39.07 & 20 & 22 \\
\hline
\end{tabular}
\end{table}

\begin{figure*}[!t]
\begin{center}
$\begin{array}{ccc}
    \includegraphics[width=0.24\textwidth]{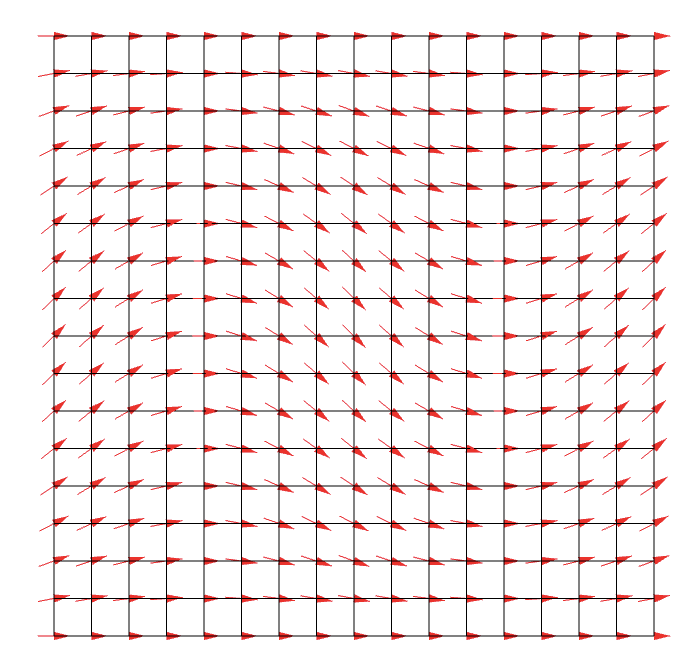} &
    \includegraphics[width=0.24\textwidth]{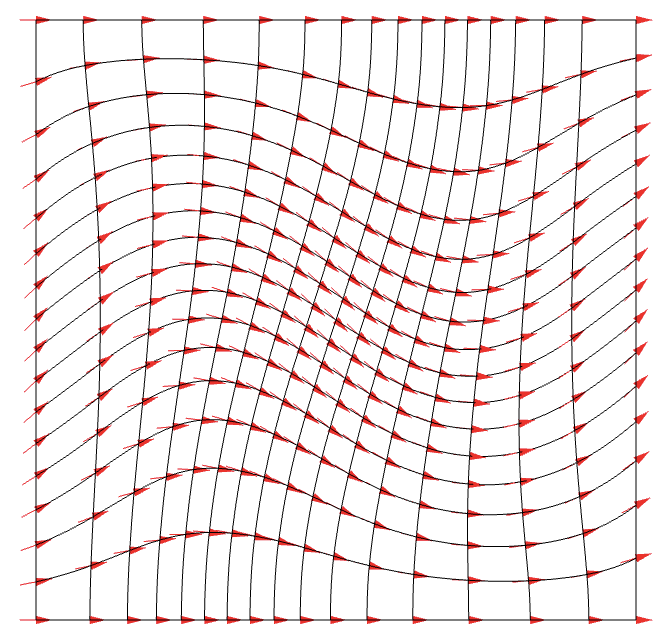} &
    \includegraphics[width=0.24\textwidth]{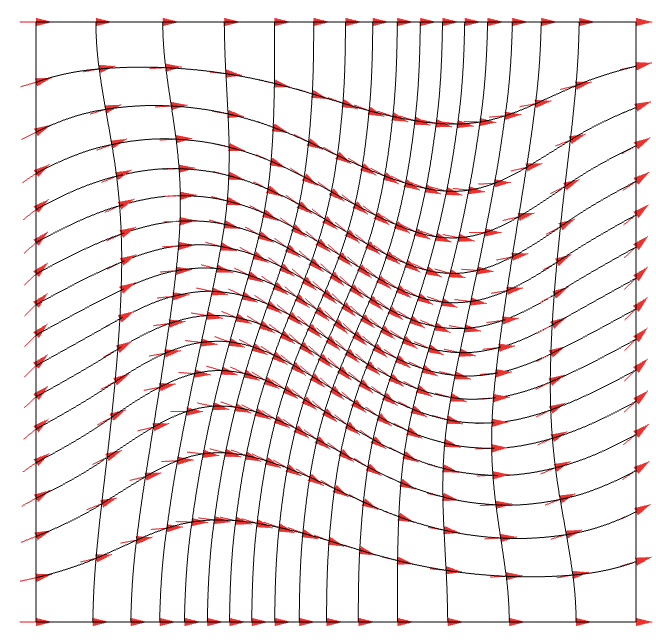} \\
    \textrm{(a)} & \textrm{(b)} & \textrm{(c)}
\end{array}$
\end{center}
\vspace{-5mm}
\caption{Shape, size and orientation adaptivity with a spatially varying target.
(a) initial mesh, (b) optimized mesh with target derivatives omitted, and
(c) optimized mesh with target derivatives included in both the
gradient and Hessian. The red arrows indicate the target orientation.}
\label{fig:discrete-target-meshes}
\end{figure*}

For $p=1$ and $p=2$, both implementations take the same numbers of
Newton iterations on CPUs and GPUs. For $p=3$, the Newton iteration
counts vary slightly as numerical round-off errors add up and lead to different convergence histories.

On CPUs, classic TMOP has lower solve times for all
three orders, benefiting from its highly optimized implementation
for the constant targets used in this benchmark.
On GPUs, the \dfem\ solve is slower for $p=1$, but faster for $p$=2 and 3.
The $p=1$ run is slower because of relatively little computational work
to leverage GPU acceleration, and this becomes apparent with speed-up for higher polynomial orders. The GPU results are obtained using the
automatically generated pointwise derivatives and the existing
\dfem\ operators, without developing TMOP-specific GPU kernels.
Note that the \dfem\ framework is under active development, and further
performance improvements are expected as its implementation is
optimized.

\subsection{\ra\ with spatially varying target $W(\vx)$}
\label{sec:variable-target-results}

This example demonstrates the convergence benefits of including the derivatives of a spatially varying target matrix during \ra. We consider the target matrix,
\begin{equation}
\label{eq:size-orientation-target}
\begin{aligned}
W(\vx)&=\alpha(\vx)R\big(\theta(\vx)\big),\\
R(\theta)&=\begin{pmatrix}
\cos\theta & -\sin\theta\\
\sin\theta & \cos\theta
\end{pmatrix},\\
\alpha(\vx)&=0.14-0.08\sin(\pi x)\sin(\pi y),\\
\theta(\vx)&=\pi y(1-y)\cos(2\pi x),
\end{aligned}
\end{equation}
which prescribes smaller elements near the center, unity aspect-ratio, $90^{\circ}$ skewness, and a varying orientation. To match the prescribed target, a shape, size, and alignment metric, $\mu_{36,vos}(T) = \frac{\lVert T-I\rVert_F^2}{\det T}$ is chosen.
The initial mesh is a uniform mesh with $16 \times 16$ quadrilateral elements, and is shown with the target orientation in Figure \ref{fig:discrete-target-meshes}(a).

\begin{table}[!bt]
\centering
\caption{\ra\ for spatially varying targets. Newton iterations $N$,
final objective $F_\mu$, and mean elapsed time $t$ (seconds) are compared for
classic TMOP and the \dfem-based approach without and with spatial derivatives of the target matrix.}
\label{tab:variable-target-timings}
\vspace{2pt}
\setlength{\tabcolsep}{3pt}
\renewcommand{\arraystretch}{1.12}
\begin{tabular}{|c|c|c|c|c|}
\hline
Method & $W$ deriv. & $N$ & $F_\mu$ & $t$ (s) \\
\hline
\multicolumn{5}{|c|}{\textbf{Analytic target}} \\
\hline
Classic & No & 17 & 2.52 & 0.38 \\
\hline
\dfem\ & No & 17 & 2.52 & 0.31 \\
\hline
\dfem\ & Yes & 6 & 2.25 & 0.21 \\
\hline
\multicolumn{5}{|c|}{\textbf{Discrete target}} \\
\hline
Classic & No & 16 & 2.60 & 0.40 \\
\hline
\dfem\ & No & 16 & 2.60 & 0.35 \\
\hline
\dfem\ & Yes & 8 & 2.36 & 0.29 \\
\hline
\end{tabular}
\end{table}

To illustrate the effect of including target derivatives, we perform experiments with both the analytic and discrete representations of the target.
For analytic targets, the target derivatives are computed directly from the analytic expressions.
For discrete targets, the target area and orientation are represented as grid functions on the initial mesh. In classic TMOP, the derivatives of the target with respect to the spatial coordinates are omitted from the gradient and the Hessian. With \dfem, these derivatives are automatically included in both the gradient and Hessian.

Table \ref{tab:variable-target-timings} reports measurements using four MPI ranks on an Apple M3 Max CPU. For both analytic and discrete targets, omitting target derivatives gives classic TMOP and \dfem\ identical convergence histories and final meshes, with lower mean solve times for \dfem. Including target derivatives allows \dfem\ to reach lower objective values in fewer Newton iterations, further reducing the solve time. These derivatives account for how node motion changes the sampled size and orientation fields, producing stronger mesh concentration near the center of the domain, as shown for the discrete target in Figure \ref{fig:discrete-target-meshes}(b,c).

\subsection{Node limiting for ALE mesh relaxation}
\label{sec:ale-node-limiting}

\begin{figure}[!t]
\begin{center}
\setlength{\arraycolsep}{1.5pt}
$\begin{array}{ccc}
    \includegraphics[width=0.32\linewidth]{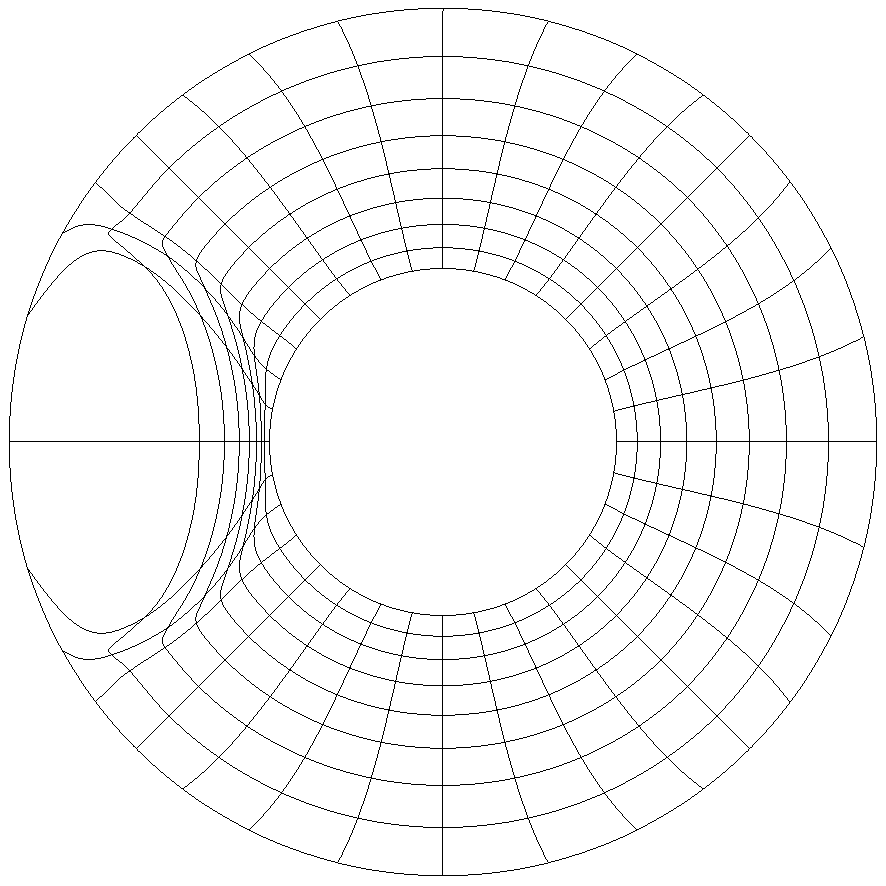} &
    \includegraphics[width=0.32\linewidth]{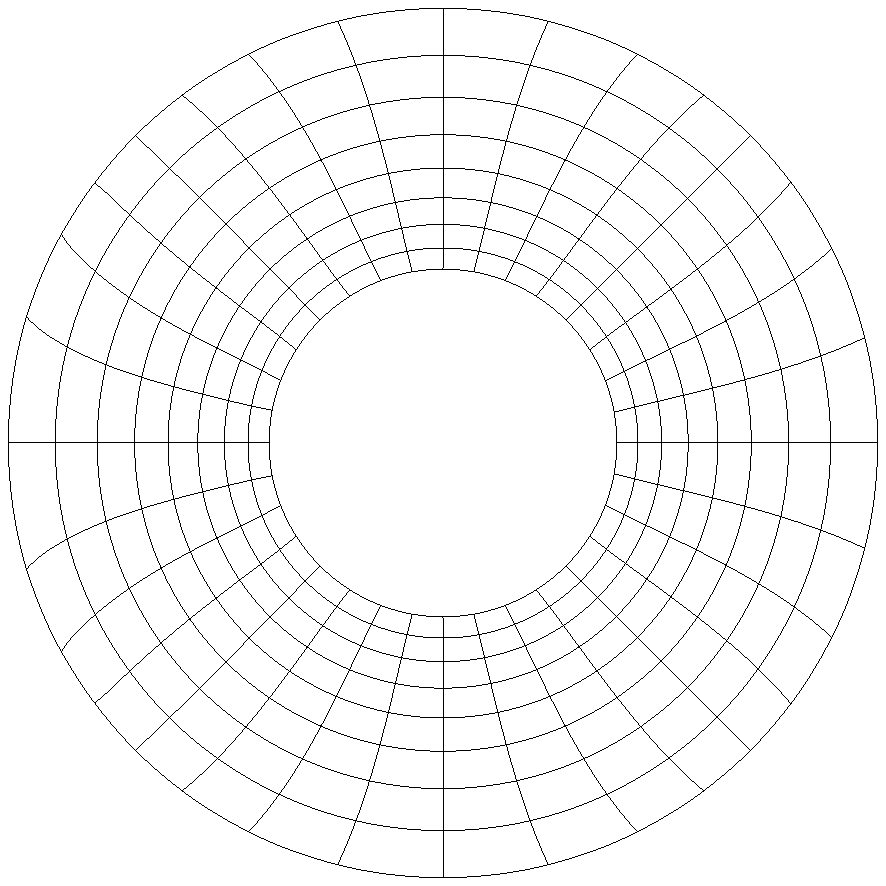} &
    \includegraphics[width=0.32\linewidth]{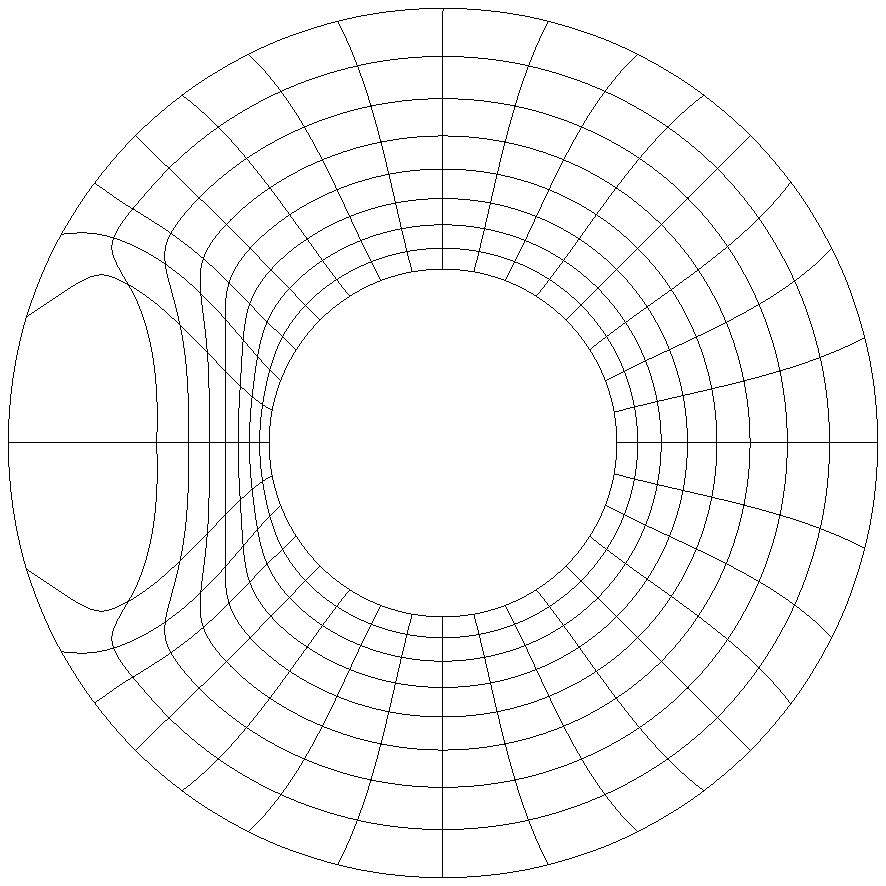} \\
    \textrm{(a)} & \textrm{(b)} & \textrm{(c)}
\end{array}$
\end{center}
\vspace{-5mm}
\caption{Node limiting during ALE mesh relaxation with \dfem\:
(a) initial Lagrangian mesh; (b) optimized mesh without limiting;
(c) optimized mesh with $c_{\mathrm{lim}}=0.5$.
Node limiting retains more of the initial Lagrangian mesh distribution
while still improving element shape.
All boundary nodes are fixed during optimization.}
\label{fig:ale-node-limiting}
\end{figure}

This example demonstrates the use of \dfem\ for ALE mesh relaxation with a node-limiting term that balances mesh-quality improvement against node displacement. Limiting node movement helps reduce numerical diffusion during the subsequent remap \cite{TMOP2020CAF}. We consider an annular cubic mesh deformed during the Lagrangian phase,  shown in Figure \ref{fig:ale-node-limiting}(a).

\begin{figure}[!t]
\centerline{\resizebox{0.8\linewidth}{!}{\begingroup\relax

\newcommand{\alePenaltyLabel}{$c_{\mathrm{lim}}$}\relax
\newcommand{\aleDisplacementLabel}{$d_{\max}$}\relax
\newcommand{\aleQualityLabel}{$F_\mu$}\relax
\newcommand{\aleOptimizedLabel}{Optimized with limiting}\relax
\newcommand{\aleUnlimitedLabel}{Optimized without limiting}\relax
\newcommand{\aleInitialLabel}{Initial mesh}\relax
\newcommand{\aleDisplacementPanelLabel}{(a)}\relax
\newcommand{\aleQualityPanelLabel}{(b)}\relax

\newcommand{\aleUnlimitedDisplacement}{0.3299341920419521}\relax
\newcommand{\aleUnlimitedQuality}{0.3023279741022578}\relax
\pgfplotstableread{
parameter displacement quality
0.001 3.151838404709964e-05 0.9998103462036682
0.01 0.0005202261681998138 0.9927310654108471
0.03 0.001848957830217099 0.9648953515565509
0.05 0.003721217197947721 0.933411283754773
0.1 0.009884699755272192 0.8541977441774292
0.2 0.02933649050566869 0.7246514520738526
0.3 0.05230918377033066 0.6334013932454041
0.5 0.0968319897019424 0.5184658469875759
1 0.1842823342292078 0.3920498584987207
2 0.270258358372735 0.3236114425064596
3 0.3003503344393527 0.3089566608036495
5 0.3186920735178712 0.3034984254148168
10 0.3270607855920588 0.3024135108021119
}\alePenaltyData\relax

\definecolor{aleDisplacementColor}{HTML}{1F5D91}\relax
\definecolor{aleQualityColor}{HTML}{1F5D91}\relax
\definecolor{aleInitialColor}{HTML}{A16D37}\relax

\begin{tikzpicture}[font=\small]
\begin{groupplot}[
  group style={group size=1 by 2, vertical sep=1.35cm},
  scale only axis,
  width=\dimexpr\linewidth-1.65cm\relax,
  height=3.4cm,
  xmode=log,
  xmin=0.001, xmax=15,
  xtick={0.001,0.01,0.1,1,10},
  axis lines=left,
  tick align=outside,
  tick label style={font=\footnotesize},
  grid=major,
  major grid style={draw=black!12},
  legend style={draw=none, fill=none, font=\scriptsize},
  legend cell align=left,
  clip=false
]

\nextgroupplot[
  ymode=log,
  ymin=0.00002, ymax=0.42,
  ytick={0.0001,0.001,0.01,0.1},
  ylabel={\aleDisplacementLabel},
  legend style={at={(0.97,0.04)}, anchor=south east}
]
\addplot[color=black!60, dashed, line width=0.8pt, mark=none]
  coordinates {(0.001,\aleUnlimitedDisplacement)
               (15,\aleUnlimitedDisplacement)};
\addlegendentry{\aleUnlimitedLabel}
\addplot[color=aleDisplacementColor, line width=1pt,
         mark=*, mark size=2pt]
  table[x=parameter,y=displacement]{\alePenaltyData};
\addlegendentry{\aleOptimizedLabel}

\nextgroupplot[
  xlabel={\alePenaltyLabel},
  ymin=0.27, ymax=1.06,
  ytick={0.3,0.4,0.5,0.6,0.7,0.8,0.9,1},
  ylabel={\aleQualityLabel},
  legend style={at={(0.03,0.06)}, anchor=south west}
]

\addplot[color=aleInitialColor, solid, line width=1pt, mark=none]
  coordinates {(0.001,1) (15,1)};
\addlegendentry{\aleInitialLabel}
\addplot[color=black!60, dashed, line width=0.8pt, mark=none]
  coordinates {(0.001,\aleUnlimitedQuality)
               (15,\aleUnlimitedQuality)};
\addlegendentry{\aleUnlimitedLabel}
\addplot[color=aleQualityColor, line width=1pt,
         mark=*, mark size=2pt]
  table[x=parameter,y=quality]{\alePenaltyData};
\addlegendentry{\aleOptimizedLabel}

\end{groupplot}
\node[anchor=north] at ([yshift=-20pt]group c1r1.south)
  {\aleDisplacementPanelLabel};
\coordinate (alePanelLabelRow) at ([yshift=-6pt]current bounding box.south);
\node[anchor=north] at (group c1r2.south |- alePanelLabelRow)
  {\aleQualityPanelLabel};
\end{tikzpicture}\relax
\endgroup\relax}}
\vspace{-5mm}
\caption{Effect of $c_{\mathrm{lim}}$ on (a) maximum nodal displacement ($d_{\max}$) and (b) normalized TMOP quality $F_\mu$. Solid lines show \ra\ with limiting and dashed lines show the corresponding values without limiting. The brown horizontal line in (b) marks the initial value $F_\mu=1$.}
\label{fig:ale-penalty-sweep}
\end{figure}
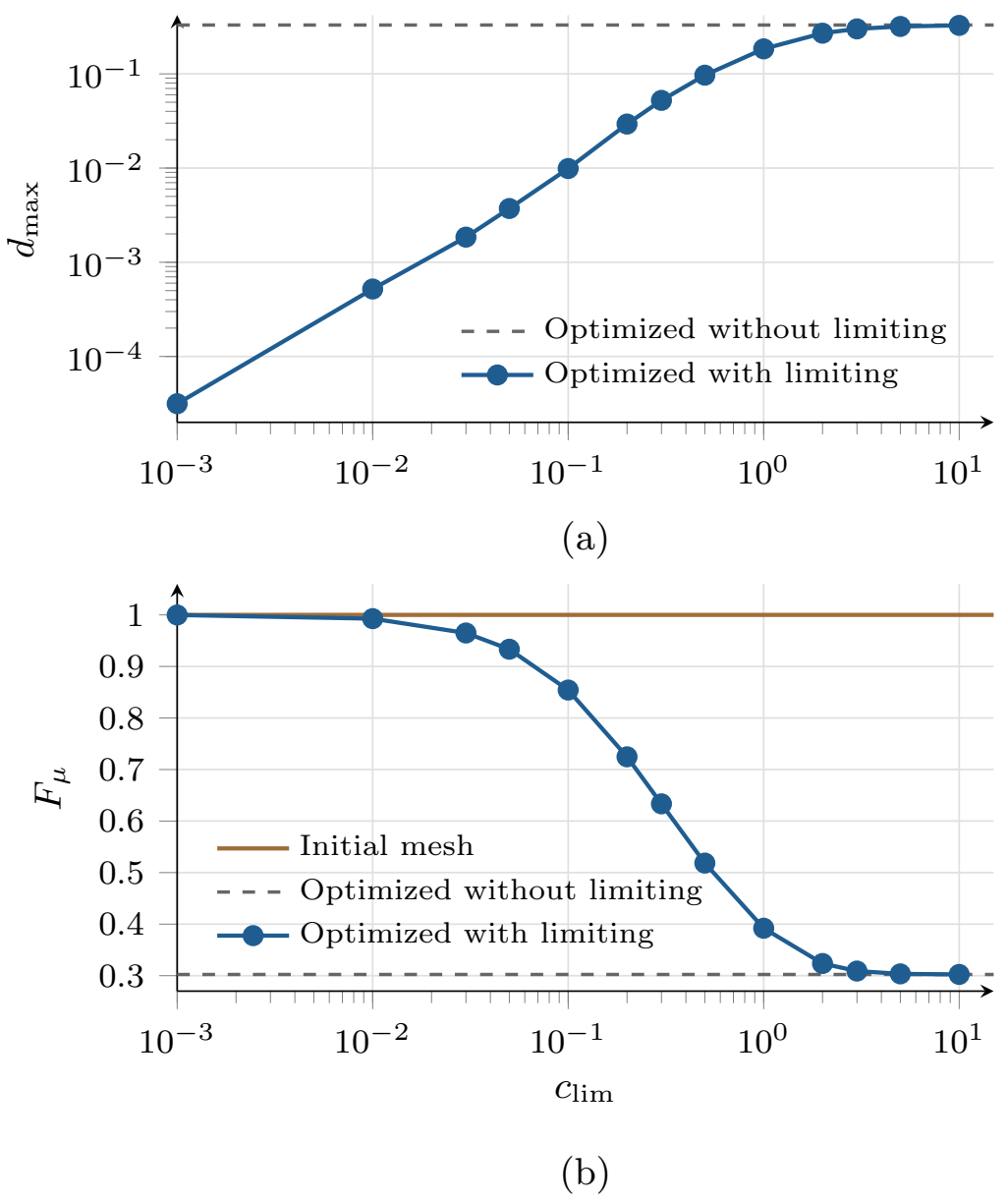

We minimize $F_\mu+F_{\mathrm{lim}}$ using the shape metric $\mu_{2,s}$, the constant ideal target $W=I$ prescribing unity aspect ratio and $90^{\circ}$ skewness, and the initial Lagrangian mesh as the reference configuration ($\vx_0$) for the displacement penalty in \eqref{eq:node-limiting}.
All boundary nodes are held fixed during \ra.
We normalize the metric contribution by its initial value, keeping this factor fixed during optimization. Thus, the initial value of $F_\mu$ is one for the results reported here. We next optimize the mesh without limiting and with limiting ($c_{\mathrm{lim}}=0.5$). The resulting meshes are shown in Figure \ref{fig:ale-node-limiting} (b,c).

Without limiting, \ra\ reduces $F_\mu$ from $1.00$ to $0.30$, accompanied by substantial redistribution of the interior nodes. With $c_{\mathrm{lim}}=0.5$, the final mesh-quality objective is $0.52$, and maximum nodal displacement ($d_{\max}$) decreases from $0.33$ without limiting to $0.10$, a reduction of $70.7\%$.
Figure \ref{fig:ale-penalty-sweep} shows how maximum nodal displacement and normalized mesh quality measure $F_\mu$ vary as $c_{\mathrm{lim}}$ is varied between $10^{-3}$ and $10$.
Decreasing $c_{\mathrm{lim}}$ strengthens node limiting and reduces $d_{\max}$ monotonically, while $F_\mu$ increases toward its initial value as less node movement is allowed.

\begin{figure*}[!tbh]
\begin{center}
\setlength{\arraycolsep}{0.3pt}
$\begin{array}{ccccc}
    \includegraphics[height=0.18\textwidth]{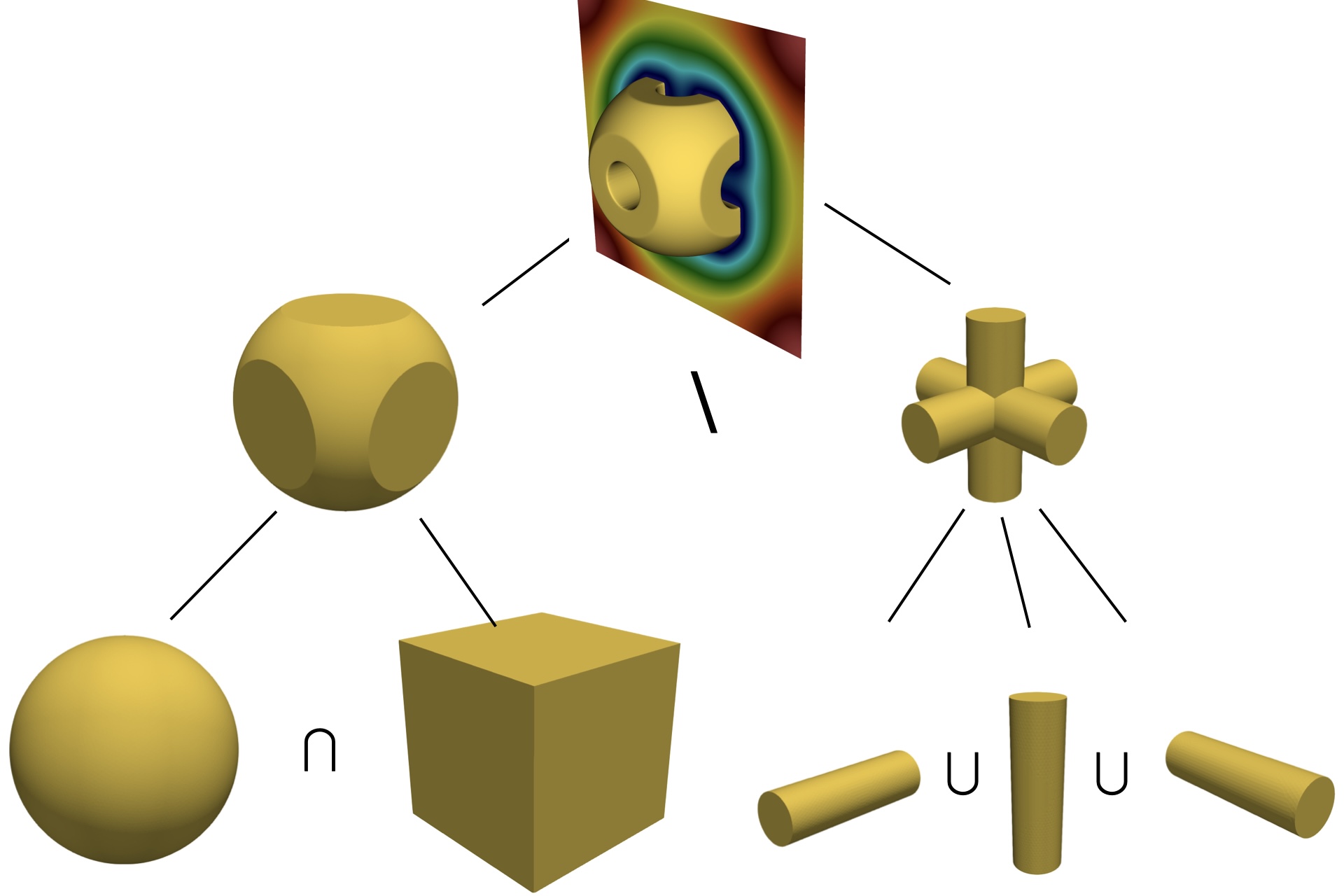} &
    \includegraphics[height=0.18\textwidth]{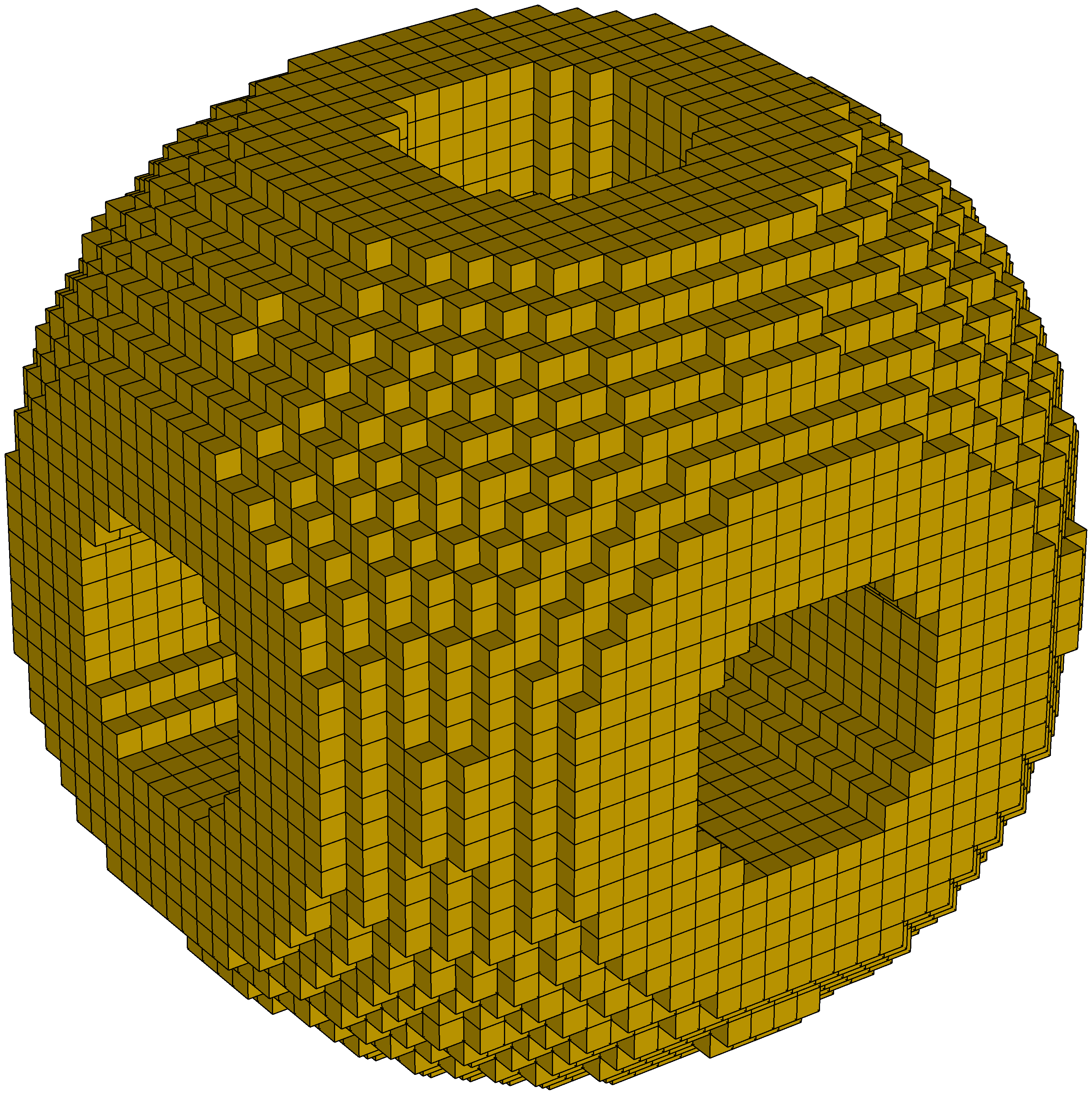} &
    \includegraphics[height=0.18\textwidth]{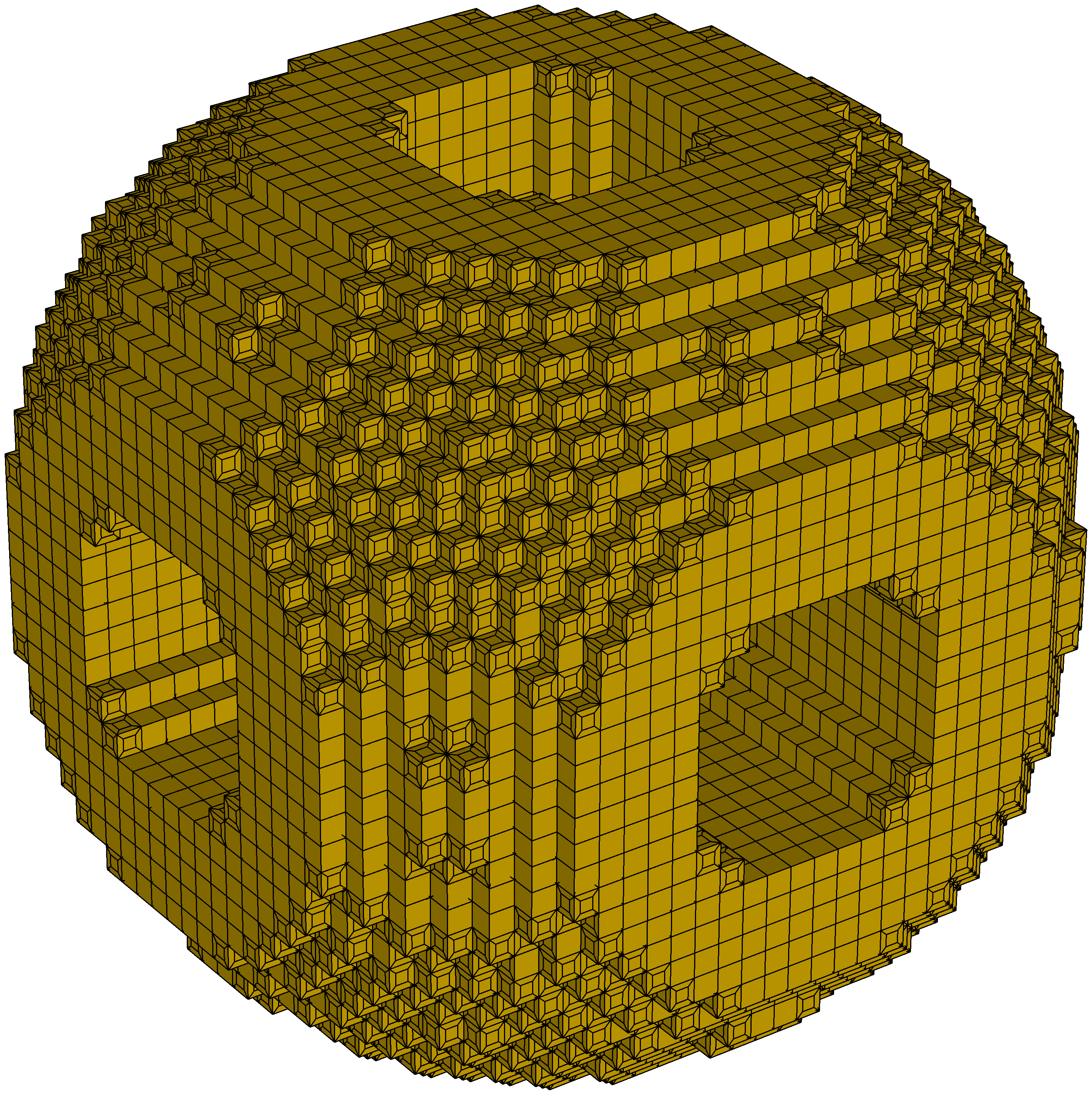} &
    \includegraphics[height=0.18\textwidth]{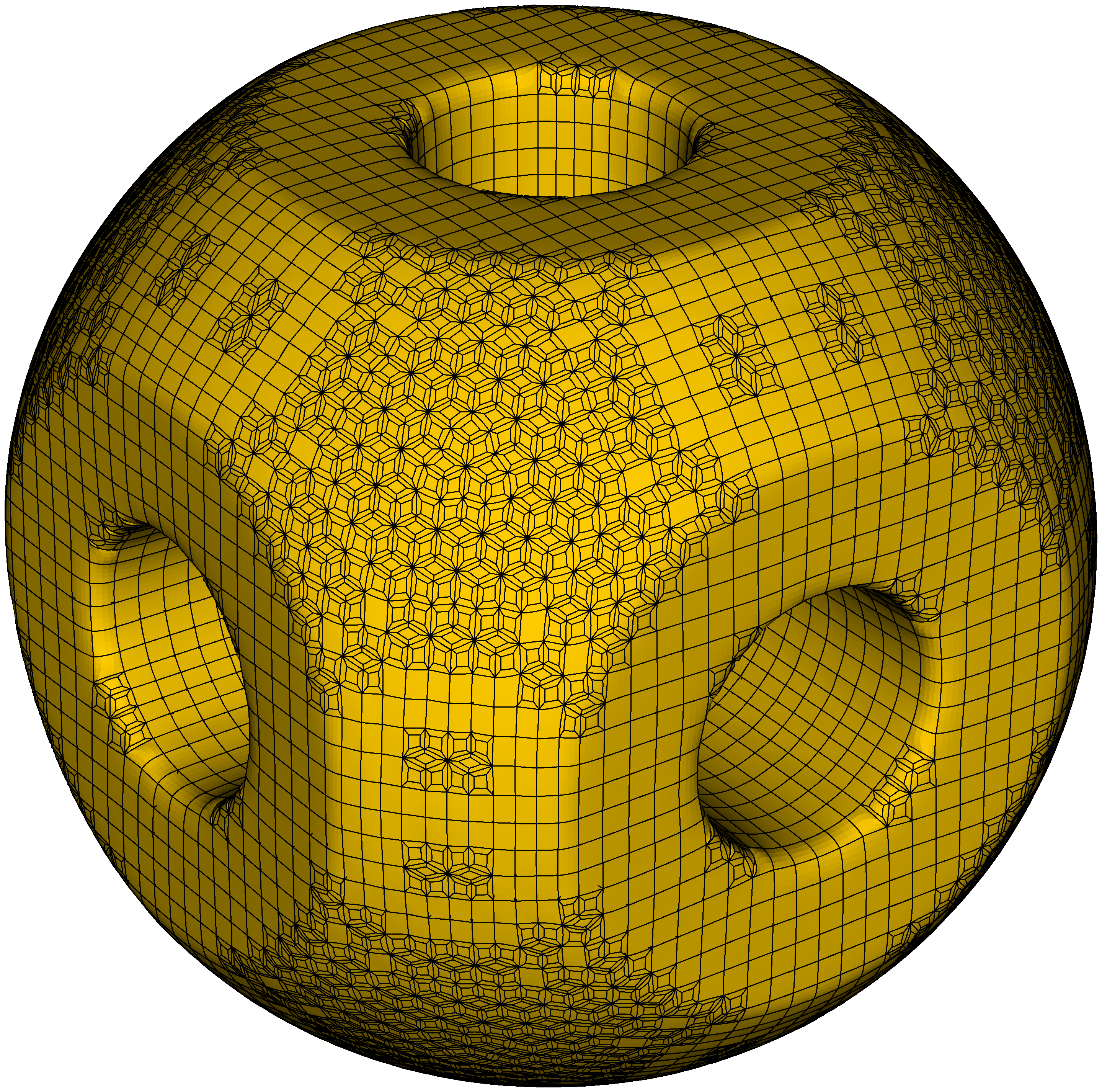} &
    \includegraphics[height=0.18\textwidth]{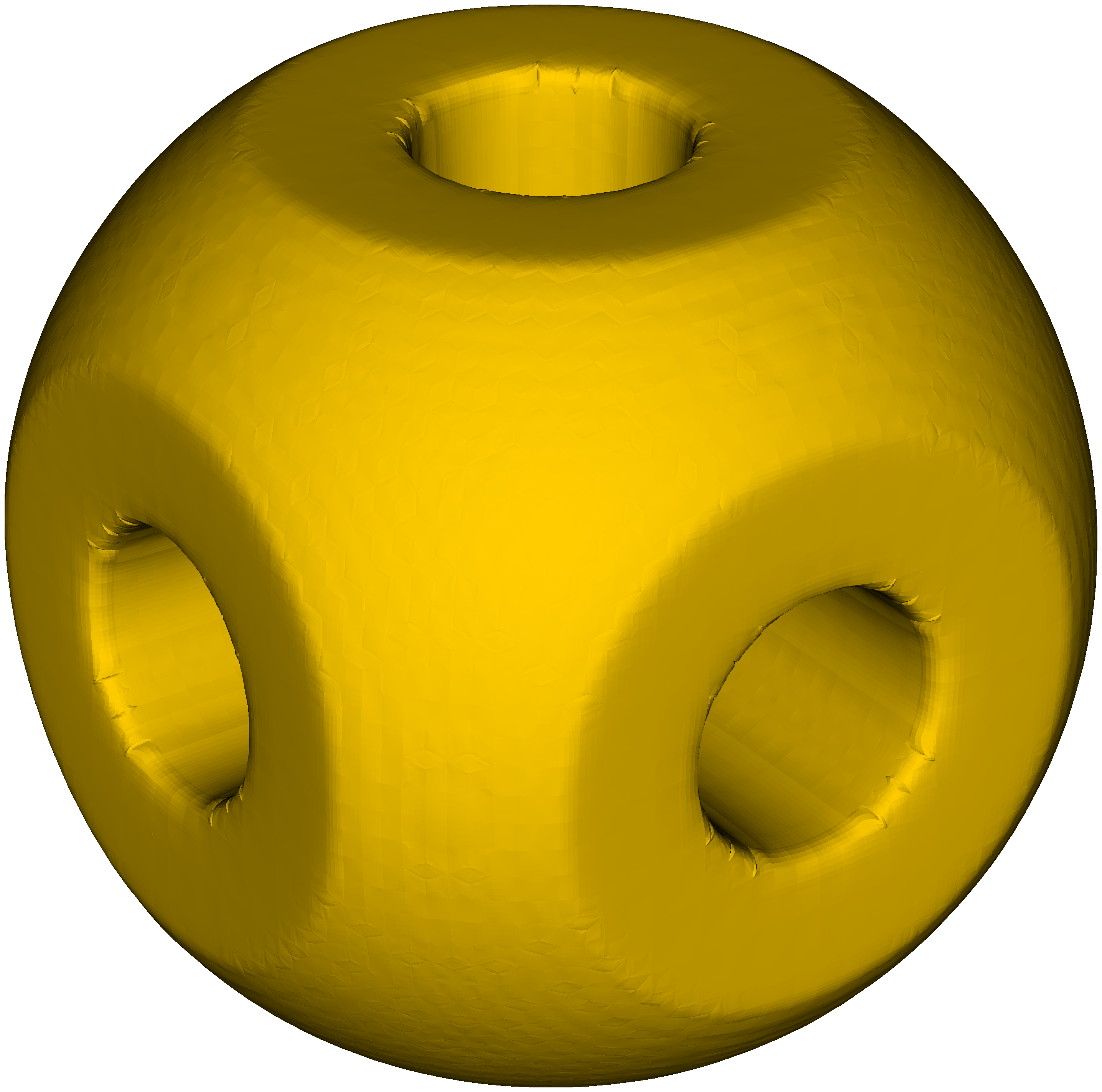} \\
    \textrm{(a)} & \textrm{(b)} & \textrm{(c)} & \textrm{(d)} & \textrm{(e)}
\end{array}$
\end{center}
\vspace{-5mm}
\caption{Boundary fitting to a CSG domain with \dfem.
(a) CSG tree defining the target geometry \cite{TMOP2023CAD}, (b) initial trimmed mesh, (c) preprocessed mesh \cite{mittal2023high}, (d) optimized mesh, and (e) its surface.}
\label{fig:csg-fitting}
\end{figure*}

\subsection{Mesh fitting for implicit geometries}
\label{sec:fitting-example}

This section shows the use of \dfem-based \ra\ to align a hexahedral mesh boundary with an implicitly defined geometry. The target geometry is represented using constructive solid geometry in Figure \ref{fig:csg-fitting}(a); a cube is intersected with a sphere, followed by removal of three mutually orthogonal cylinders. The level-set function capturing the resulting boundary is also shown in Figure \ref{fig:csg-fitting}(a) and the reader is referred to Section 4.2 of \cite{TMOP2023CAD} for details of its construction.

Starting with a Cartesian-aligned quadratic hexahedral mesh, elements outside the domain are trimmed to obtain a stair-stepped mesh, as shown in Figure \ref{fig:csg-fitting}(b). This mesh is preprocessed following \cite{mittal2023high}
to refine hex elements at the boundary such that \ra\ for interface fitting will not flatten their faces or edges. The preprocessed mesh is shown in Figure \ref{fig:csg-fitting}(c).
Next, the TMOP functional with level-set-based penalization ($F_\mu+F_{\mathrm{surf}}$) is minimized. For the mesh fitting term, $\mathcal S$ contains the boundary nodes and $\sigma$ represents the level-set function.
The mesh-quality term is set up to use a target based on the Jacobian of the initial mesh, and a shape and size quality metric, $\mu(T)=\lVert T-T^{-\mathsf T}\rVert_F^2$.
The optimized mesh in Figure \ref{fig:csg-fitting}(d) conforms to the
curved outer boundary and cylindrical openings of the target geometry,
which is evident with a pointwise error of $10^{-6}$ for the final mesh.
The surface of the optimized mesh is shown in Figure \ref{fig:csg-fitting}(e).

Note that classic TMOP does not support partial assembly and GPU acceleration for the mesh fitting functional $F_{\mathrm{surf}}$ as those specialized kernels have not been developed yet. In contrast, the proposed \dfem-based framework readily supports GPU acceleration for evaluating $F_{\mathrm{surf}}$ and its derivatives.

\subsection{Sensitivity analysis of the \ra\ operator}
\label{sec:digital-twin}

This section considers a simple example of solving the heat conduction problem in a multimaterial domain with a parameterized interface, to demonstrate that the \dfem-based framework can provide consistent sensitivities associated with the \ra\ operator.

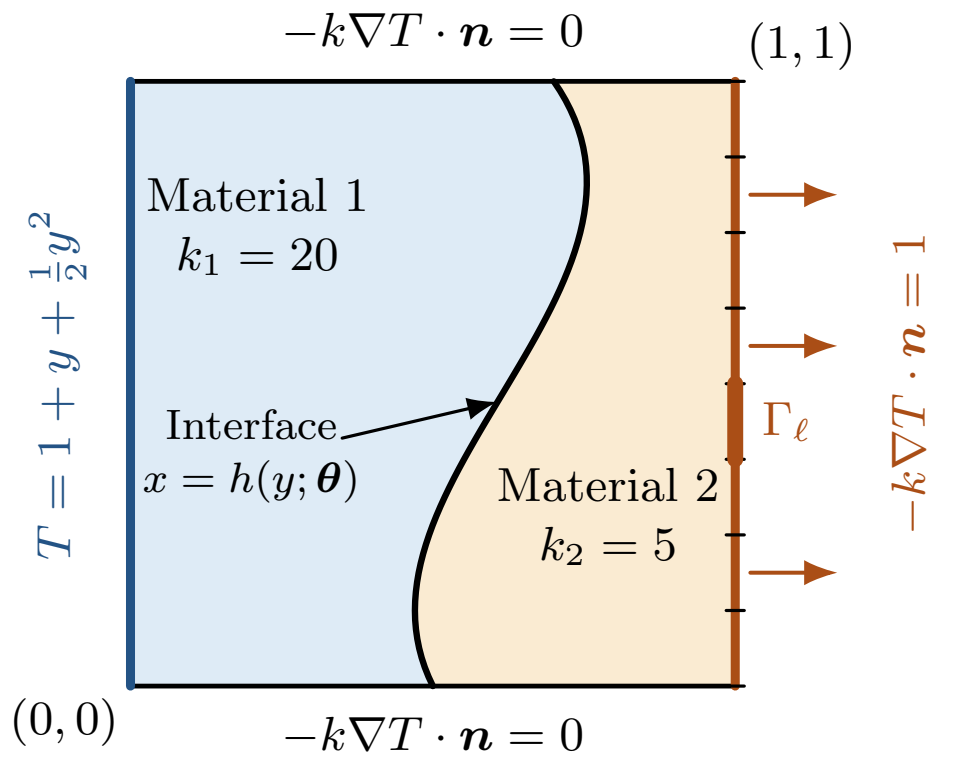
\begin{figure}[!t]
\centerline{\begingroup
\definecolor{twinleft}{RGB}{222,235,246}
\definecolor{twinright}{RGB}{250,235,210}
\definecolor{twintemperature}{RGB}{35,83,133}
\definecolor{twinflux}{RGB}{170,78,22}
\begin{tikzpicture}[
  x=4.0cm, y=4.0cm, font=\small,
  line cap=round, line join=round,
  >=Latex,
  declare function={interface(\t)=0.6+0.2*(\t-0.5)+0.1*\t*(1-\t)
                                  +1.6*(\t-0.5)*\t*(1-\t);}
]
  \fill[twinleft] (0,0) -- (0,1) -- ({interface(1)},1)
    -- plot[domain=1:0, samples=65, variable=\t]
      ({interface(\t)},\t) -- cycle;
  \fill[twinright] ({interface(0)},0)
    -- plot[domain=0:1, samples=65, variable=\t]
      ({interface(\t)},\t) -- (1,1) -- (1,0) -- cycle;

  \draw[line width=0.9pt] (0,0) rectangle (1,1);
  \draw[twintemperature, line width=1.7pt] (0,0) -- (0,1);
  \draw[line width=1.2pt]
    plot[domain=0:1, samples=65, variable=\t]
      ({interface(\t)},\t);

  \node[align=center] at (0.21,0.76) {Material 1\\$k_1=20$};
  \node[align=center] at (0.79,0.28) {Material 2\\$k_2=5$};
  \node[align=center, font=\footnotesize] at (0.20,0.38)
    {Interface\\$x=h(y;\boldsymbol\theta)$};
  \draw[->, line width=0.6pt] (0.35,0.41)
    -- ({interface(0.47)},0.47);

  \draw[twinflux, line width=1.7pt] (1,0) -- (1,1);
  \foreach \j in {0,...,8}
    {\draw[line width=0.6pt] (0.985,{\j/8}) -- (1.015,{\j/8});}
  \draw[twinflux, line width=3pt] (1,0.375) -- (1,0.5);
  \node[twinflux, anchor=west, inner sep=2pt] at (1.025,0.4375)
    {$\Gamma_\ell$};
  \foreach \y in {0.1875,0.5625,0.8125}
    {\draw[twinflux, ->, line width=0.75pt] (1.025,\y) -- (1.17,\y);}

  \node[twintemperature, rotate=90] at (-0.12,0.5)
    {$T=1+y+\tfrac12 y^2$};
  \node[twinflux, rotate=90] at (1.29,0.5)
    {$-k\nabla T\cdot\boldsymbol n=1$};
  \node[align=center, anchor=south] at (0.5,1.02)
    {$-k\nabla T\cdot\boldsymbol n=0$};
  \node[align=center, anchor=north] at (0.5,-0.02)
    {$-k\nabla T\cdot\boldsymbol n=0$};
  \node[anchor=north east, inner sep=2pt] at (0,0) {$(0,0)$};
  \node[anchor=south west, inner sep=2pt] at (1,1) {$(1,1)$};
\end{tikzpicture}
\endgroup}
\vspace{-5mm}
\caption{Heat digital twin: reference interface and boundary conditions.
Ticks mark eight sensor segments; $\boldsymbol n$ is the outward unit normal.}
\label{fig:digital-twin-domain}
\end{figure}

\begin{figure*}[!tb]
\begin{center}
\setlength{\arraycolsep}{1pt}
$\begin{array}{ccccc}
    \includegraphics[width=0.23\textwidth]{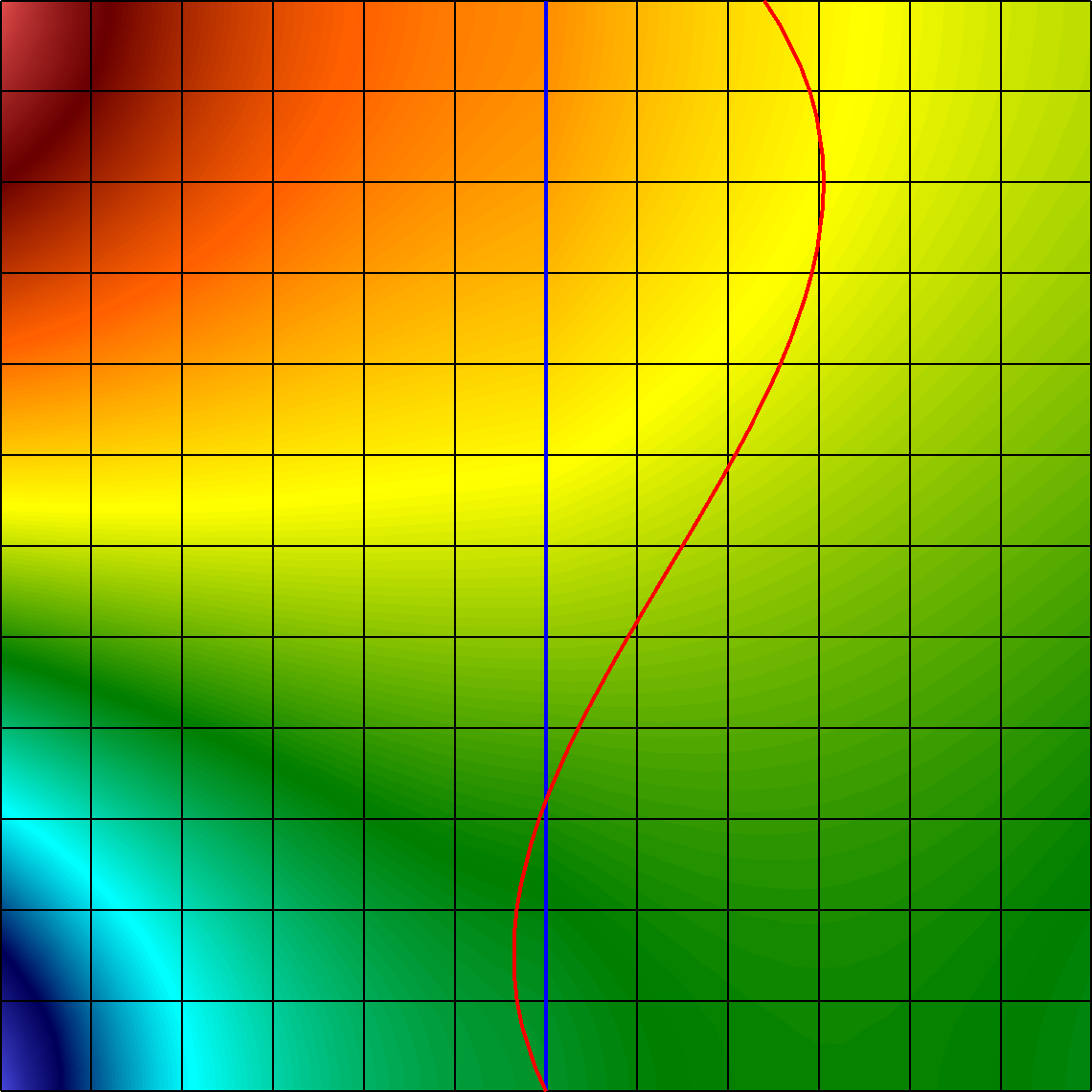} &
    \includegraphics[width=0.23\textwidth]{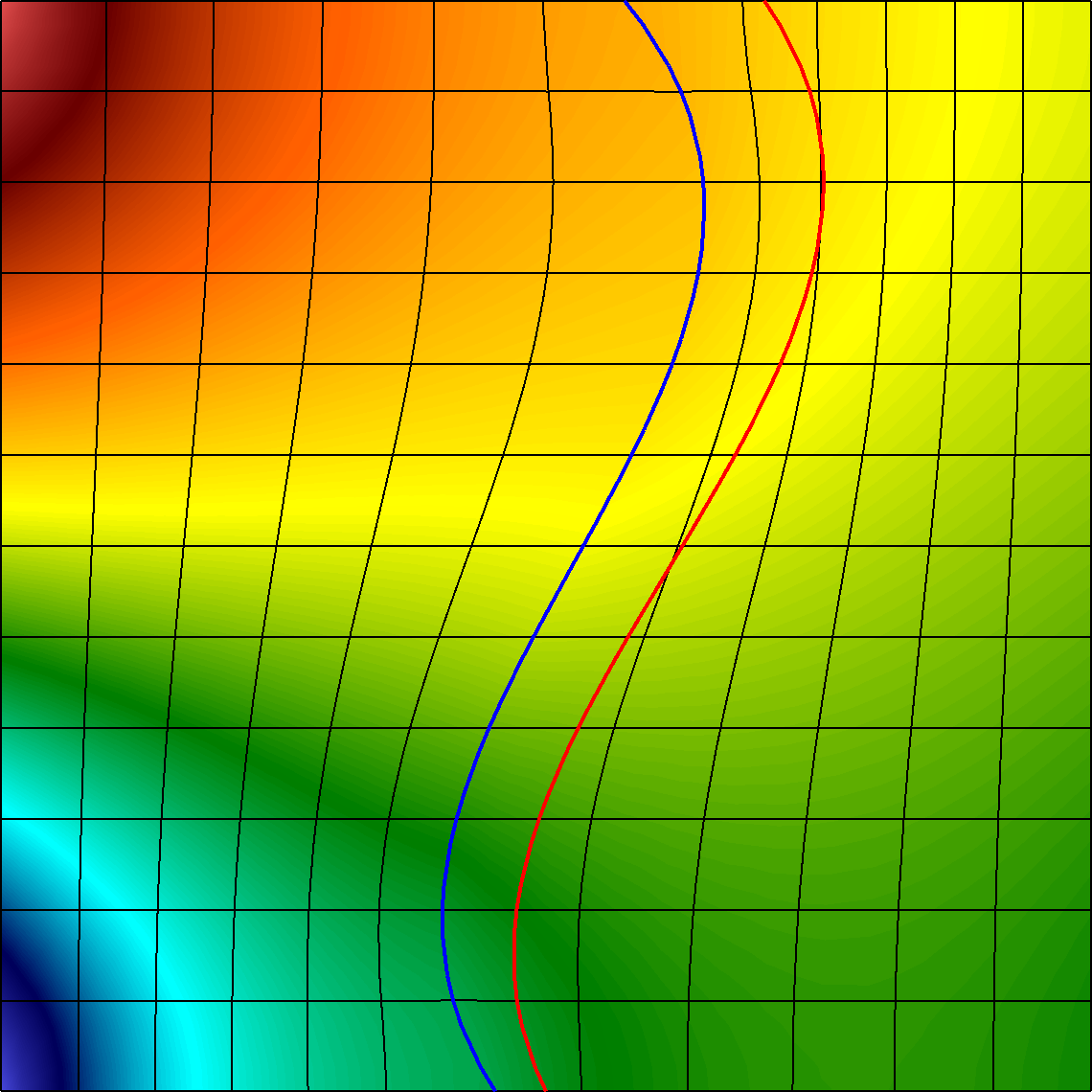} &
    \includegraphics[width=0.23\textwidth]{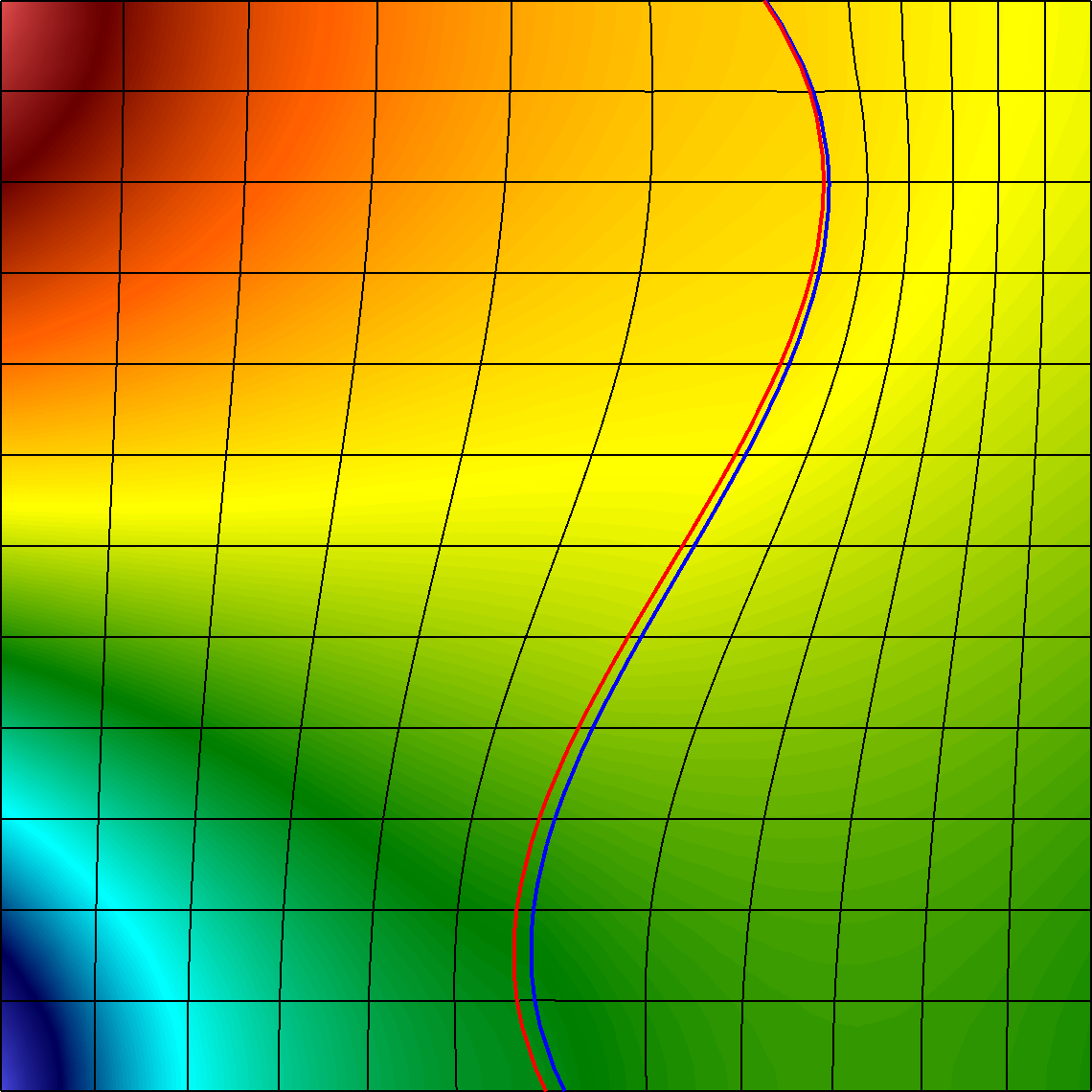} &
    \includegraphics[width=0.23\textwidth]{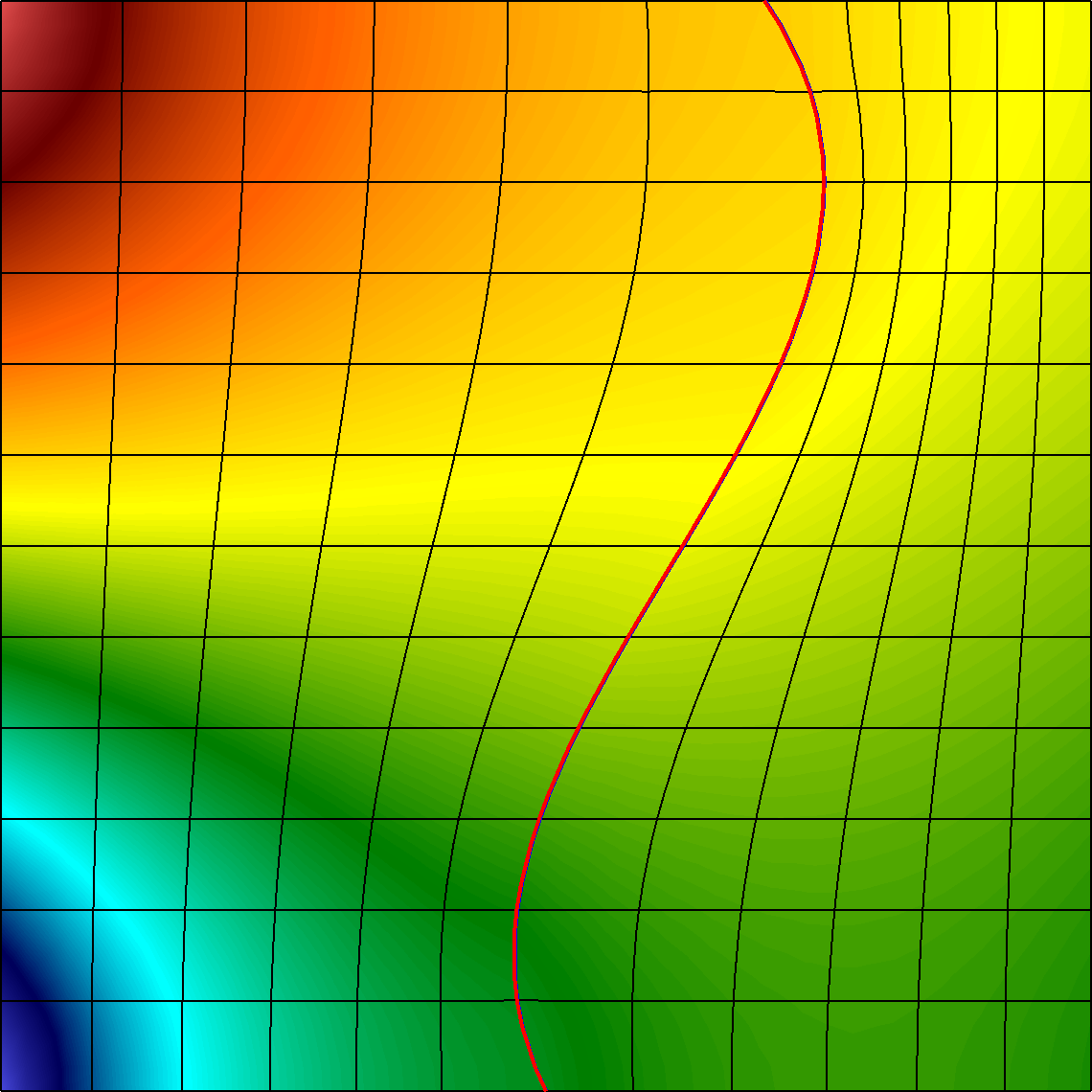} &
    \includegraphics[width=0.05\textwidth]{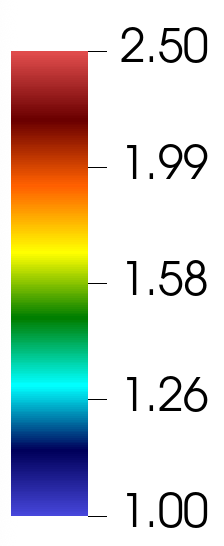}\\
    \textrm{(a)}\,\,n=0 & \textrm{(b)}\,\,n=1 &
    \textrm{(c)}\,\,n=2 & \textrm{(d)}\,\,n=3 &
\end{array}$
\end{center}
\vspace{-5mm}
\caption{Temperature fields and fitted meshes during interface recovery at
inverse-design (Gauss-Newton) iteration $n$; $n=0$ denotes the initial guess.
Red and blue curves denote the reference and current interfaces,
respectively. Nodes on $x=0$ and $x=1$ are held fixed.}
\label{fig:digital-twin-recovery}
\end{figure*}

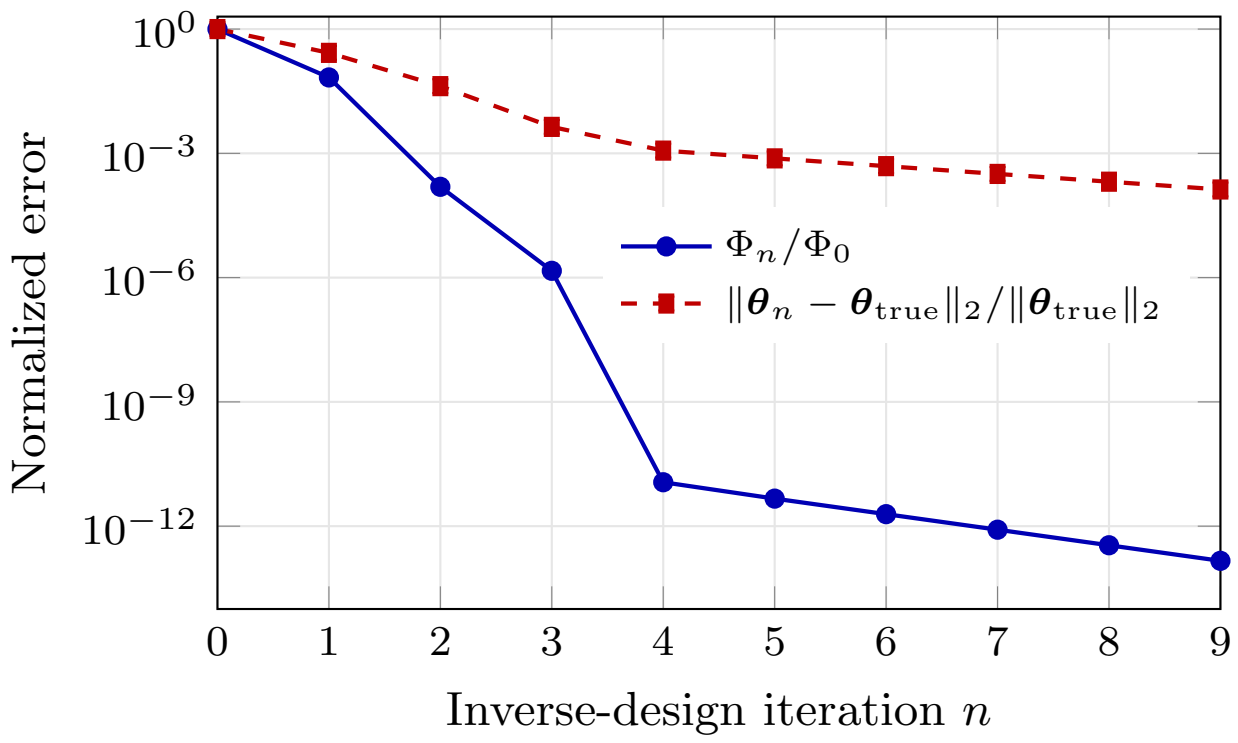
\begin{figure}[!btp]
\centerline{\begin{tikzpicture}
\begin{semilogyaxis}[
  width=0.97\linewidth, height=5.5cm,
  xlabel={Inverse-design iteration $n$}, ylabel={Normalized error},
  xmin=0, xmax=9, ymin=1e-14, ymax=2,
  xtick={0,1,2,3,4,5,6,7,8,9},
  ytick={1e0,1e-3,1e-6,1e-9,1e-12},
  grid=major, grid style={gray!20},
  tick label style={font=\footnotesize},
  label style={font=\small},
  legend style={font=\scriptsize, draw=none, fill=white,
    at={(0.97,0.68)}, anchor=north east, cells={anchor=west}},
  mark size=1.5pt,
]
\addplot[blue!70!black, thick, mark=*]
  table[x=iteration,y=relative_objective,col sep=comma]
  {figures/digital-twin.csv};
\addlegendentry{$\Phi_n/\Phi_0$}
\addplot[red!75!black, thick, dashed, mark=square*]
  table[x=iteration,y=relative_parameter_error,col sep=comma]
  {figures/digital-twin.csv};
\addlegendentry{$\|\vtheta_n-\vtheta_{\mathrm{true}}\|_2/\|\vtheta_{\mathrm{true}}\|_2$}
\end{semilogyaxis}
\end{tikzpicture}}
\vspace{-5mm}
\caption{Normalized sensor misfit and relative parameter error versus
inverse-design (Gauss-Newton) iteration $n$.}
\label{fig:digital-twin-convergence}
\end{figure}

Figure \ref{fig:digital-twin-domain} shows the computational domain, $\Omega=[0,1]^2$, with a known Dirichlet boundary condition ($T=1+y+y^2/2$) on the left boundary and a homogeneous Neumann condition on the top and bottom boundaries.
The multimaterial interface $\sigma(\vx;\vtheta)=x-h(y;\vtheta)=0$ is parameterized as
\begin{equation}
\label{eq:digital-twin-interface}
\begin{aligned}
h(y;\vtheta)={}&\tfrac12+\theta_0+\theta_1(y-\tfrac12)
                 +\theta_2y(1-y)\\
              &+4\theta_3(y-\tfrac12)y(1-y),
\end{aligned}
\end{equation}
such that $\vtheta$ controls its deformation. The two materials have conductivities of 20 and 5.

The idea with this demonstration is to first start with some initial value of $\vtheta$ and then minimize $F_{\mu} + F_{\mathrm{surf}}$ to obtain an interface-fitted mesh. Next, the heat conduction problem ($-\nabla\cdot(k\nabla T)=0$) is solved and the resulting temperature ($T$) on the right boundary is compared with a provided \emph{target}. The target specifies average temperatures on $L$ equal-length panels covering the right boundary. The computed panel averages are $s_\ell=|\Gamma_\ell|^{-1}\int_{\Gamma_\ell}T\,\dd s$. The overarching goal is to determine $\vtheta$ that minimizes the objective:
\begin{equation}
\label{eq:digital-twin-objective}
\Phi(\vtheta)=\frac12\sum_{\ell=1}^{L}
\big(s_\ell(\vtheta)-s_\ell^{\mathrm{obs}}\big)^2,
\end{equation}
which entails consistent sensitivity analysis of the \ra\ operator and the PDE solve.

For this numerical experiment, we set $L=8$ and use a $12\times12$ quadrilateral mesh. The initial mesh with material interface at $x=0.5$, i.e. $\vtheta=0$, is shown in Figure \ref{fig:digital-twin-recovery}(a). Here, we also show the curved interface that is the ground truth resulting from $\vtheta_{\mathrm{true}} = (0.1,0.2,0.1,0.4)$.
Note that the shown mesh connectivity and material distribution are held fixed irrespective of $\vtheta$ so as to separate the sensitivity of the \ra\ operator from non-smooth operations such as topology or material changes. Additionally, both the mesh nodal coordinates and temperature are represented using $p=2$ bases. To simplify the sensitivity calculation, the mesh nodes on $x=0$ and $x=1$ are held fixed. Nodes on $y=0$ and $y=1$ may move tangentially to accommodate the material-interface endpoints.

Minimizing the sensor objective using the Gauss-Newton method requires
the sensor sensitivities and the objective gradient with respect to
$\vtheta$. Let $T$ denote the discrete temperature solution,
$\boldsymbol s$ the sensor values, and
$\mathcal R(T,\vX^*)=\boldsymbol0$ denote the discrete heat
equations. The sensor values depend on both the temperature solution
and the fitted mesh, so the chain rule gives
\begin{equation}
\label{eq:digital-twin-sensitivity-chain}
\begin{aligned}
\frac{\dd\Phi}{\dd\vtheta}
={}&\left(\frac{\dd\boldsymbol s}{\dd\vtheta}\right)^{\mathsf T}
\big(\boldsymbol s-\boldsymbol s^{\mathrm{obs}}\big),\\
\frac{\dd\boldsymbol s}{\dd\vtheta}
={}&\frac{\partial\boldsymbol s}{\partial T}
\frac{\dd T}{\dd\vtheta}
+\frac{\partial\boldsymbol s}{\partial\vX}
\frac{\dd\vX^*}{\dd\vtheta}.
\end{aligned}
\end{equation}
The fixed left and right boundary nodes have zero coordinate sensitivities.
Consequently, the mesh-motion contributions from the prescribed boundary
data and sensor quadrature vanish, and
$\frac{\partial\boldsymbol s}{\partial\vX}
\frac{\dd\vX^*}{\dd\vtheta}=\boldsymbol0$ in this example.
The temperature and mesh sensitivities follow by differentiating the
heat equations and the TMOP stationarity condition
\eqref{eq:mesh-stationarity}, respectively:
\begin{equation}
\label{eq:digital-twin-tangent-equations}
\begin{aligned}
\frac{\partial\mathcal R}{\partial T}
\frac{\dd T}{\dd\vtheta}
+\frac{\partial\mathcal R}{\partial\vX}
\frac{\dd\vX^*}{\dd\vtheta}&=\boldsymbol0,\\
\frac{\partial\mathcal J}{\partial\vX}
\frac{\dd\vX^*}{\dd\vtheta}
+\frac{\partial\mathcal J}{\partial\vtheta}&=\boldsymbol0.
\end{aligned}
\end{equation}
Substituting these sensitivities into
\eqref{eq:digital-twin-sensitivity-chain} gives
\begin{equation}
\label{eq:digital-twin-sensitivity-expanded}
\begin{aligned}
\frac{\dd\boldsymbol s}{\dd\vtheta}
={}&\left[
\frac{\partial\boldsymbol s}{\partial T}
\left(\frac{\partial\mathcal R}{\partial T}\right)^{-1}
\frac{\partial\mathcal R}{\partial\vX}
-\frac{\partial\boldsymbol s}{\partial\vX}
\right]
\left(\frac{\partial\mathcal J}{\partial\vX}\right)^{-1}
\frac{\partial\mathcal J}{\partial\vtheta}.
\end{aligned}
\end{equation}
We use \dfem\ to compute the heat-residual mesh derivative and the TMOP
Hessian and mixed derivative in \eqref{eq:digital-twin-sensitivity-expanded}.

Figure \ref{fig:digital-twin-recovery} shows how the predicted interface is virtually indistinguishable from the true interface in just three Gauss-Newton updates of $\vtheta$ as $\Phi$ decreases from $7.36\times10^{-3}$ to $1.07\times10^{-8}$.
Figure \ref{fig:digital-twin-convergence} shows the convergence of the sensor objective ($\Phi_n/\Phi_0$) and relative error in $\vtheta$ with Gauss-Newton updates. At termination after nine updates, $\Phi=1.07\times10^{-15}$ and the relative parameter error is $1.34\times10^{-4}$ ($0.013\%$).
The sensitivity of the objective to the parameters considered here is rather low, due to the nature of the problem setup, but the example nonetheless shows the effectiveness of the \dfem-based approach for consistent sensitivity analysis of the \ra\ operator.

\section{Conclusions}
\label{sec:conclusions}
This work presents a framework for differentiable \ra\ of high-order meshes. By leveraging \dfem, this framework requires us to simply specify the pointwise mesh-adaptation objective. These local functions are automatically differentiated with \dfem, and combined into a global gradient vector and a partially assembled Hessian operator. The generated functional evaluation and derivative kernels readily support CPU and GPU execution.
Numerical experiments show GPU solve times
comparable to those of hand-derived TMOP implementations on a benchmark problem, improvement in $r$-adaptivity convergence due to accurate computation of the derivatives, and demonstrate the use of mesh sensitivities for material-interface identification in an inverse heat-conduction problem.

Future work will integrate the proposed framework in ALE simulations and PDE-constrained optimization. We will continue to leverage the latest advances in \dfem\ to further improve the performance of \ra. Finally, we will also look at ways to compute consistent sensitivities in problems where the mesh topology changes.

\section*{Acknowledgments}
\anon{
lorem ipsum dolor sit amet consectetur adipiscing elit nisi dolor mollit dolores officia labore in nisi eum quod et incididunt et sunt dolore consequat velit officia laborum autem excepturi consectetur maxime et dolore dolor quas soluta in consectetur expedita facere dolor proident ducimus culpa do accusamus excepteur repellendus minim quis et quidem laboris ullamco in enim exercitation in aliquip nobis est dolorum elit lorem velit illum libero adipiscing cillum in sed et labore consequat non quidem minus id lorem corrupti
}
{
This manuscript has been authored by Lawrence Livermore National Security, LLC under Contract No. DE-AC52-07NA27344 with the U.S. Department of Energy and was supported by the LLNL-LDRD Program under project No 26-SI-002. The United States Government retains, and the publisher, by accepting the article for publication, acknowledges that the United States Government retains a non-exclusive, paid-up, irrevocable, world-wide license to publish or reproduce the published form of this manuscript, or allow others to do so, for United States Government purposes. LLNL-CONF-2024904.
}

\clearpage

\bibliographystyle{ieeetr}
\bibliography{imr}

\end{document}